\documentclass[12pt,reqno]{article}
\usepackage{amssymb}
\usepackage{amsthm}
\usepackage{amsmath}
\newtheorem{theorem}{Theorem}[section]

\newtheorem{definition}[theorem]{Definition}
\newtheorem{example}[theorem]{Example}

\begin{document}
\begin{center} \Large{\bf On invariants of second-order ordinary differential equations $y''=f(x,y,y')$ via point transformations}
\end{center}
\medskip
\begin{center}
Ahmad Y. Al-Dweik\\
Department of Mathematics \& Statistics, King Fahd University of Petroleum and Minerals, Dhahran 31261, Saudi Arabia\\
aydweik@kfupm.edu.sa
\end{center}
\begin{abstract}
Bagderina \cite{Bagderina2013} solved the equivalence problem for
a family of scalar second-order ordinary differential equations
(ODEs), with cubic nonlinearity in the first-order derivative, via
point transformations. However, the question is open for the
general class $y''=f(x,y,y')$ which is not cubic in the
first-order derivative. We utilize Lie's infinitesimal method to
study the differential invariants of this general class under an
arbitrary point equivalence transformations. All fifth order
differential invariants and the invariant differentiation
operators are determined. As an application, invariant description
of all the canonical forms in the complex plane for second-order
ODEs $y''=f(x,y,y')$ where both of the two Tress\'e relative
invariants are non-zero is provided.
\end{abstract}
\section{Introduction}
Lie's group classification of ODEs shows that the second-order
equations can possess one, two, three or eight infinitesimal
symmetries. According to Lie's classification \cite{Lie1888} in
the complex domain, any second order ODE
\begin{equation}\label{w1}
\begin{array}{ll}
y''=f(x,y,y'),\\
\end{array}
\end{equation}
is obtained by a change of variables from one of eight canonical
forms. The equations with eight symmetries and only these
equations can be linearized by a change of variables. The initial
seminal studies of scalar second-order ODEs which are linearizable
by means of point transformations are due to Lie \cite{Lie1888}
and Tress\'e \cite{Tresse1894}. They showed that the latter
equations are at most cubic in the first derivative and gave a
convenient invariant description of all linearizable equations.
Mahomed and Leach \cite{Mahomed1989} proved that Lie linearization
conditions are equivalent to the vanishing of the Tress\'e
relative invariants (\ref{wer}) as stated in the next theroem
\begin{theorem}\cite{Mahomed2007}
The equation $ y'' = f(x,y,y')$ is equivalent to the normal form
$y'' = 0$ with eight symmetries under \textbf{point
transformations} if and only if the Tress\'e relative invariants
\begin{equation}\label{wer}
\begin{array}{ll}
I_1=f_{y',y',y',y'}\\
I_2=\dot{D}_x^2 f_{y',y'}-4 \dot{D}_x f_{y,y'}-3 f_{y}f_{y',y'} +6f_{y,y}+f_{y'} \left( 4f_{y,y'}  - \dot{D}_x f_{y',y'} \right) \\
\end{array}
\end{equation}
both vanish identically. Where $ \dot{D}_x = \frac{\partial }
{{\partial x}} + y'\frac{\partial }{{\partial y}}+f\frac{\partial}
{{\partial y'}}$.
\end{theorem}
Regarding the equivalence of  the second-order differential
equations to the normal form $y'' = 0$ via contact
transformations, it is well known that all second-order
differential equations have $y'' = 0$ as the sole equivalence
class.

The linearization problem is a particular case of the equivalence
problem. For the general theory of the  equivalence problem
including algorithms of construction of differential invariants,
the interested reader is referred to \cite{Olver1995,Ovsiannikov}.
Ibragimov \cite{Ibr1997,Ibr1999} developed a simple method for
constructing invariants of families of linear and nonlinear
differential equations admitting infinite equivalence
transformation groups. Lie's infinitesimal method was applied to
solve the equivalence problem for several linear and nonlinear
differential equations
\cite{Ibr2002-1,Ibr2002-2,Ibr2007,Ibr2004,Ibr2005,Johnpillai2001,Torrisi2004,Torrisi2005,Tracina2004}.
Cartan's equivalence method \cite{Olver1995,Gardner1989} is
another systematic approach to solve the equivalence problem for
differential equations.

By using Lie's infinitesimal method, Bagderina
\cite{Bagderina2013} solved the equivalence problem of
second-order ODEs which are at most cubic  in the first-order
derivative ($I_1=0$)
\begin{equation}
\begin{array}{ll}
y''=a(x,y){y'}^3+b(x,y){y'}^2+c(x,y)y'+d(x,y)\\
\end{array}
\end{equation}
with respect to the group of point equivalence transformations
\begin{equation}\label{o2}
\begin{array}{ll}
\bar{x}=\phi(x,y), \bar{y}=\psi(x,y).\\
\end{array}
\end{equation}
As an extension, in this paper, we use Lie's infinitesimal method
to study the differential invariants of the second-order ODEs
(\ref{w1}) which are not cubic in the first-order derivative
($I_1\ne0$) with respect to the group of point equivalence
transformations. The motivation of this study is finding invariant
description of the canonical forms for second-order ODEs in the
complex plane \cite{Olver1995} which are not cubic in the
first-order derivative.

Invariant description of the canonical forms for second-order ODEs
in the complex plane with three infinitesimal symmetries was given
in \cite{Ibr2008, Ibr2007} where they presented the candidates for
all four types and then they studied these candidates. In this
paper, invariant description of all the canonical forms in the
complex plane for second-order ODEs $y''=f(x,y,y')$ where both of
the two Tress\'e relative invariants (\ref{wer}) are non-zero is
provided.

The structure of the paper is as follows. In the next section, we
give a short description of Lie's infinitesimal method to find the
differential invariants and invariant differentiation operators of
the class of ODEs (\ref{w1}) with respect to the group of point
equivalence transformations (\ref{o2}). In Section 3, using the
methods described in Section 2, the infinitesimal point
equivalence transformations are recovered. In Section 4, we find
the fifth-order differential invariants and invariant
differentiation operators of the class of ODEs (\ref{w1}), which
are not cubic in the first-order derivative, under point
equivalence transformations. In Section 5, invariant description
of all the canonical forms in the complex plane \cite{Olver1995}
for second-order ODEs $y''=f(x,y,y')$ where both of the two
Tress\'e relative invariants (\ref{wer}) are non-zero is provided.
Finally, the conclusion is presented.

Throughout this paper, we use the notation $A=[a_1,a_2,...,a_{n}]$
to express any differential operator
$A=\displaystyle{\sum^{n}_{j=1} a_j \frac{\partial}{\partial
b_j}}$. Also, we denote $y'$ by $p$.
\section{Lie's infinitesimal method}
In this section, we  briefly describe the Lie method used
to derive differential invariants using point equivalence
transformations.

Consider the $k$th-order system of
PDEs of $n$ independent variables
 $x = (x^{1}, x^{2}, ... , x^{n})$ and $m$ dependent variables $u = (u^{1}, u^{2},...,
 u^{m})$
 \begin{equation}\label{a1}
 \begin{array}{ll}
 E_{\alpha}(x,u,...,u_{(k)})=0, & \alpha=1,...,m~,
\end{array}
\end{equation}
where $u_{(1)},u_{(2)},...,u_{(k)}$ denote the collections of all
first, second, ..., $k$th-order partial derivatives, i.e.,
$u_{i}^{\alpha}=D_{i}(u^{\alpha}),\;
u_{ij}^{\alpha}=D_{j}D_{i}(u^{\alpha})$,..., respectively, in
which the total differentiation operator with respect to $x^{i}$
is
\begin{equation}\label{a2}
\begin{array}{ll}
D_{i}=\frac{\partial}{\partial
x^{i}}+u_{i}^{\alpha}\frac{\partial}{\partial
u^\alpha}+u_{ij}^\alpha\frac{\partial}{\partial
u_{j}^\alpha}+...,&i=1,...,n~,
\end{array}
\end{equation}
with the summation convention adopted for repeated indices.
\begin{definition}\rm
\emph{The Lie-B\"acklund operator} is
\begin{equation}\label{a4}
\begin{array}{ll}
X=\xi^{i}\frac{\partial}{\partial
{x^i}}+\eta^\alpha\frac{\partial}{\partial
u^\alpha}&\xi^i,\eta^\alpha\in \emph{A}~,
\end{array}
\end{equation}
where $A$ is the space of  \emph{differential functions}.

The operator (\ref{a4}) is an abbreviated form of the infinite formal
sum
\begin{equation}\label{a5}
\begin{array}{ll}
 X^{(s)}&=\xi^{i}\frac{\partial}{\partial x^{i}}+\eta^{\alpha}\frac{\partial}{\partial u^{\alpha}}+\sum\limits_{s\geq 1 }\zeta_{i_{1}i_{2}...i_{s}}^{\alpha}\frac{\partial}{\partial u_{i_1 i_2...i_s}^\alpha},\\
  &=\xi^{i}D_i+W^\alpha\frac{\partial}{\partial u^{\alpha}}+\sum\limits_{s\geq 1 }D_{i_1}...D_{i_s}(W^\alpha)\frac{\partial}{\partial u_{i_1 i_2...i_s}^\alpha},\\
\end{array}
\end{equation}
where the additional coefficients are determined uniquely by the
prolongation formulae
\begin{equation}\label{a6}
\begin{array}{ll}
\zeta_{i}^\alpha=D_{i}(\eta^{\alpha})-u_{j}^\alpha D_{i}(\xi^{j})=D_{i}(W^\alpha)+\xi^{j}u_{ij}^\alpha,\\
\zeta_{i_1...i_s}^\alpha=D_{i_s}(\zeta_{i_1...i_{s-1}}^\alpha)-u_{ji_1...i_{s-1}}^\alpha D_{i_s}(\xi^{j})=D_{i_1}...D_{i_s}(W^\alpha)+\xi^j u_{ji_1...i_s}^\alpha,& s>1,\\
\end{array}
\end{equation}
in which $W^\alpha$ is the\emph{ Lie characteristic function}
\begin{equation}\label{a7}
W^\alpha=\eta^\alpha-\xi^j u_{j}^{\alpha}.
\end{equation}
\end{definition}
\begin{definition}\rm
\emph{The point equivalence transformation of a class of PDEs
(\ref{a1})} is an invertible transformation of the independent and
dependent variables of the form
\begin{equation}\label{a7}
\begin{array}{ll}
\bar{x}=\phi(x,u), \bar{u}=\psi(x,u),\\
\end{array}
\end{equation}
that maps every equation of the class into an equation of the
same family, viz.
\begin{equation}\label{a8}
\begin{array}{ll}
E_{\alpha}(\bar{x},\bar{u},...,\bar{u}_{(k)})=0, & \alpha=1,...,m.
\end{array}
\end{equation}
\end{definition}
In order to describe Lie's infinitesimal method for deriving
differential invariants using point equivalence transformations,
we use as example the class of equations (\ref{w1}). It is well-known that the
point equivalence transformation
\begin{equation}\label{a9}
\begin{array}{ll}
\bar{x}=\phi(x,y), \bar{y}=\psi(x,y),\\
\end{array}
\end{equation}
 maps (\ref{w1}) into the same family, viz.
\begin{equation}\label{a10}
\bar{y}''=\bar{f}(\bar{x},\bar{y},\bar{y}'),
\end{equation}
for arbitrary functions  $\phi(x,y)$ and $\psi(x,y)$, where
$\bar{f}$, in general, can be different from the original function
$f$. The set of all equivalence transformations forms a group
denoted by $\mathcal{E}$.

The standard procedure for Lie's infinitesimal invariance
criterion \cite{Ovsiannikov} is implemented in the next section to
recover the continuous group of point equivalence transformations
(\ref {a9}) for the class of second-order ODEs (\ref{w1}) with the
corresponding infinitesimal point equivalence transformation
operator
\begin{equation}\label{a11}
Y=\xi(x,y)D_x+W\partial_y+D_x(W)\partial_{p}+\mu(x,y,p,f)\partial_f,\\
\end{equation}
where  $\xi(x,y)$ and $\eta(x,y)$ are arbitrary functions obtained from
\begin{equation}\label{a12}
\bar{x}= x+\epsilon~\xi(x,y)+O(\epsilon^2)=\phi(x,y),\\
\end{equation}
\begin{equation}\label{a13}
\bar{y}= y+\epsilon~\eta(x,y)+O(\epsilon^2)=\psi(x,y) ,\\
\end{equation}
and
\begin{equation}\label{a14}
\mu=\dot{D}_x^2(W)+\xi(x,y)\dot{D}_x f,\\
\end{equation} with
$W=\eta-\xi p$ and $ \dot{D}_x = \frac{\partial } {{\partial x}} +
p\frac{\partial }{{\partial y}}+f\frac{\partial} {{\partial p}}$.
\begin{definition}\rm
\emph{An invariant of a class of second-order ODEs (\ref{w1})} is
a function of the form
 \begin{equation}\label{a15}
\begin{array}{ll}
J=J(x,y,p,f),
\end{array}
\end{equation}
which is invariant under the equivalence transformation
(\ref{a9}).
\end{definition}
\begin{definition}\rm
\emph{A differential invariant of order $s$ of a class of
second-order ODEs (\ref{w1})} is a function of the form
 \begin{equation}\label{a16}
\begin{array}{ll}
J=J(x,y,p,f,f_{(1)},f_{(2)},...,f_{(s)}),\\
\end{array}
\end{equation}
which is invariant under the equivalence transformation (\ref{a9})
where $f_{(1)},f_{(2)},...,f_{(s)}$ denote the collections of all
first, second,..., $s$th-order partial derivatives.
\end{definition}
\begin{definition}\rm
\emph{An invariant system of order $s$ of a class of second-order
ODEs (\ref{w1})} is the system of the form $E_{\alpha}(x,y,p,f,f_{(1)},f_{(2)},...,f_{(s)})=0, ~ \alpha=1,...,m~$ which satisfies the condition\\
 \begin{equation}\label{a22}
\begin{array}{ll}
Y^{(s)}E_{\alpha}(x,y,p,f,f_{(1)},f_{(2)},...,f_{(s)})=0~(mod~ E_{\alpha}=0,~\alpha=1,...,m), ~ \alpha=1,...,m.\\
\end{array}
\end{equation}
An invariant system with $\alpha=1$ is called an invariant
equation.
\end{definition}
Now, according to the theory of invariants of infinite
transformation groups \cite{Ovsiannikov}, the invariant criterion
\begin{equation}\label{a17}
Y J(x,y,p,f)=0,\\
\end{equation}
should be split by means of the functions  $\xi(x,y)$ and $\eta(x, y)$ and
their derivatives. This gives rise to a homogeneous linear system
of PDEs whose solution gives the required invariants.

It should be noted that since the generator $Y$ contains
arbitrary functions $\xi(x,y)$ and $\eta(x, y)$, the
corresponding identity (\ref{a17}) leads to $m$ linear PDEs for
$J$, where $m$ is the number of the arbitrary functions and their
derivatives that appear in $Y$. We point out that these $m$ PDEs
are not necessarily linearly independent.

In order to determine the differential invariants of order $s$, we
need to calculate the prolongations of the operator $Y$ using
(\ref{a5}) by considering $f$  as a dependent variable and the
variables $x,y,p$ as independent variables:
\begin{equation}\label{a18}
\begin{array}{ll}
Y^{(s)}&=Y(x)\tilde{D}_x+Y(y)\tilde{D}_y+Y(p)\tilde{D}_p+\tilde{W}\frac{\partial}{\partial f}+\sum\limits_{s\geq 1 }\tilde{D}_{i_1}...\tilde{D}_{i_s}(\tilde{W})\frac{\partial}{\partial f_{i_1 i_2...i_s}},\\
&i_1, i_2,..., i_s\in \{x,y,p\},\\
\end{array}
\end{equation}
where
\begin{equation}\label{a19}
\begin{array}{ll}
\tilde{D}_k= \partial_{k}+f_k \partial_{f}+f_{k i} \partial_{f_{i}}+f_{k i j} \partial_{f_{ij}}+...,& i,j,k\in \{x,y,p\}. \\
\end{array}
\end{equation}
in which $\tilde{W}$ is the\emph{ Lie characteristic function}
\begin{equation}\label{a20}
\tilde{W}=\mu-Y(x) f_{x}-Y(y) f_{y}-Y(p) f_{p}.
\end{equation}
The differential invariants are determined by the equations
\begin{equation}\label{a21}
Y^{(s)} J(x,y,p,f,f_{(1)},f_{(2)},...,f_{(s)})=0.\\
\end{equation}
It should be noted that since the generator $Y^{(s)}$
contains arbitrary functions $\xi(x,y)$ and $\eta(x, y)$,  the
corresponding identity (\ref{a21}) leads to $m$ linear PDEs for
$J$, where $m$ is the number of the arbitrary functions and their
derivatives that appear in $Y^{(s)}$.

For simplicity, from here on, we denote the derivative of
$f(x,y,p)$ with respect to the independent variables $x,y,p$
 as $f_1,f_2,f_3$. The same notation will be used for higher-order derivatives.

Now, in order to find all the fifth order differential invariants
of the third-order ODE (\ref{w1}), one can solve the invariant
criterion (\ref{a21}) with $s=5$. However, for compactness of the
derived differential invariants, one can replace any partial
derivative with respect to $x$ by the total derivative with
respect to $x$. So, we need to solve the following invariant
criterion
\begin{small}
\begin{equation}\label{c1}
\begin{array}{ll}
Y^{(5)}J(x,y_{{}},y_{{1}},f_{{}},f_{{2}},f_{{3}},f_{{2,2}},f_{{2,3}},f_{{3,3}},f_{{2,2,2}},f_{{2,2,3}},f_{{2,3,3}},f_{{3,3,3}},f_{{2,2,2,2}},f_{{2,2,2,3}},f_{{2,2,3,3}},f_{{2,3,3,3}},f_{{3,3,3,3}},\\
f_{{2,2,2,2,2}},f_{{2,2,2,2,3}},f_{{2,2,2,3,3}},f_{{2,2,3,3,3}},f_{{2,3,3,3,3}},f_{{3,3,3,3,3}},d_{{1,1}},d_{{1,2}},d_{{1,3}},d_{{1,4}},d_{{1,5}},d_{{1,6}},d_{{1,7}},d_{{1,8}},d_{{1,9}},d_{{1,10}},\\
d_{{1,11}},d_{{1,12}},d_{{1,13}},d_{{1,14}},d_{{1,15}},d_{{2,1}},d_{{2,2}},d_{{2,3}},d_{{2,4}},d_{{2,5}},d_{{2,6}},d_{{2,7}},d_{{2,8}},d_{{2,9}},d_{{2,10}},d_{{3,1}},d_{{3,2}},d_{{3,3}},d_{{3,4}},\\
d_{{3,5}},d_{{3,6}},d_{{4,1}},d_{{4,2}},d_{{4,3}},d_{{5,1}})=0\\
\end{array}
\end{equation}
\end{small}
by prolonging the infinitesimal operator $Y^{(5)}$ to the
variables $d_{i,j}$ through the infinitesimals $Y^{(5)}d_{i,j}$,
where
\begin{small}
\begin{equation}\label{c2}
\begin{array}{ll}
d_{{1,1}}=\dot{D}_x{f},d_{{1,2}}=\dot{D}_x{f_2},d_{{1,3}}=\dot{D}_x{f_3},d_{{1,4}}=\dot{D}_x{f_{2,2}},d_{{1,5}}=\dot{D}_x{f_{2,3}},d_{{1,6}}=\dot{D}_x{f_{3,3}},d_{{1,7}}=\dot{D}_x{f_{2,2,2}},\\
d_{{1,8}}=\dot{D}_x{f_{2,2,3}},d_{{1,9}}=\dot{D}_x{f_{2,3,3}},d_{{1,10}}=\dot{D}_x{f_{3,3,3}},d_{{1,11}}=\dot{D}_x{f_{2,2,2,2}},d_{{1,12}}=\dot{D}_x{f_{2,2,2,3}},d_{{1,13}}=\dot{D}_x{f_{2,2,3,3}},\\
d_{{1,14}}=\dot{D}_x{f_{2,3,3,3}},d_{{1,15}}=\dot{D}_x{f_{3,3,3,3}},d_{{2,1}}=\dot{D}_x^2{f},d_{{2,2}}=\dot{D}_x^2{f_2},d_{{2,3}}=\dot{D}_x^2{f_3},d_{{2,4}}=\dot{D}_x^2{f_{2,2}},\\
d_{{2,5}}=\dot{D}_x^2{f_{2,3}},d_{{2,6}}=\dot{D}_x^2{f_{3,3}},d_{{2,7}}=\dot{D}_x^2{f_{2,2,2}},d_{{2,8}}=\dot{D}_x^2{f_{2,2,3}},d_{{2,9}}=\dot{D}_x^2{f_{2,3,3}},d_{{2,10}}=\dot{D}_x^2{f_{3,3,3}},\\
d_{{3,1}}=\dot{D}_x^3{f},d_{{3,2}}=\dot{D}_x^3{f_2},d_{{3,3}}=\dot{D}_x^3{f_3},d_{{3,4}}=\dot{D}_x^3{f_{2,2}},d_{{3,5}}=\dot{D}_x^3{f_{2,3}},d_{{3,6}}=\dot{D}_x^3{f_{3,3}},\\
d_{{4,1}}=\dot{D}_x^4{f},d_{{4,2}}=\dot{D}_x^4{f_2},d_{{4,3}}=\dot{D}_x^4{f_3},d_{{5,1}}=\dot{D}_x^5{f}\\
\end{array}
\end{equation}
\end{small}
\begin{definition}\rm
\emph{An invariant differentiation operator of a class of
second-order ODEs (\ref{w1})} is a differential operator
$\mathcal{D}$ which satisfies that if $I$ is a differential
invariant of ODE (\ref{w1}), then $\mathcal{D} I$ is its
differential invariant too.
\end{definition}
As it is shown in \cite{Ovsiannikov}, the number of independent
invariant differentiation operators $\mathcal{D}$ equals  the
number of independent variables $x,y$ and $p$. The invariant
differentiation operators $\mathcal{D}$ should take the form
\begin{equation}
\begin{array}{ll}\label{cc}
\mathcal{D}=K\tilde{D}_x+L\tilde{D}_y+M\tilde{D}_p,\\
\end{array}
\end{equation}
with the coordinates $K, L$ and $M$ satisfying the non-homogeneous
linear system
\begin{equation}
\begin{array}{llll}\label{ccc}
Y^{(5)}{K}=\mathcal{D}(Y(x)),&Y^{(5)}{L}=\mathcal{D}(Y(y)),&Y^{(5)}{M}=\mathcal{D}(Y(p)),\\
\end{array}
\end{equation}
where $K, L$ and $M$ are functions of the following variables
\begin{small}
\begin{equation}\label{cccc}
\begin{array}{ll}
x,y_{{}},y_{{1}},f_{{}},f_{{2}},f_{{3}},f_{{2,2}},f_{{2,3}},f_{{3,3}},f_{{2,2,2}},f_{{2,2,3}},f_{{2,3,3}},f_{{3,3,3}},f_{{2,2,2,2}},f_{{2,2,2,3}},f_{{2,2,3,3}},f_{{2,3,3,3}},f_{{3,3,3,3}},\\
f_{{2,2,2,2,2}},f_{{2,2,2,2,3}},f_{{2,2,2,3,3}},f_{{2,2,3,3,3}},f_{{2,3,3,3,3}},f_{{3,3,3,3,3}},d_{{1,1}},d_{{1,2}},d_{{1,3}},d_{{1,4}},d_{{1,5}},d_{{1,6}},d_{{1,7}},d_{{1,8}},d_{{1,9}},d_{{1,10}},\\
d_{{1,11}},d_{{1,12}},d_{{1,13}},d_{{1,14}},d_{{1,15}},d_{{2,1}},d_{{2,2}},d_{{2,3}},d_{{2,4}},d_{{2,5}},d_{{2,6}},d_{{2,7}},d_{{2,8}},d_{{2,9}},d_{{2,10}},d_{{3,1}},d_{{3,2}},d_{{3,3}},d_{{3,4}},\\
d_{{3,5}},d_{{3,6}},d_{{4,1}},d_{{4,2}},d_{{4,3}},d_{{5,1}}
\end{array}
\end{equation}
\end{small}
In reality, the general solution of the system (\ref{ccc}) gives
both  the differential invariants and the differential operators.
This general solution can be found by prolonging the infinitesimal
operator $Y^{(5)}$ to the variables $K, L$ and $M$ through the
infinitesimals $Y^{(5)}{K}, Y^{(5)}{L}$ and $Y^{(5)}{M}$
respectively. Then solving the invariant criterion
\begin{footnotesize}
\begin{equation}\label{ccccc}
\begin{array}{ll}
Y^{(5)}J(x,y_{{}},y_{{1}},f_{{}},f_{{2}},f_{{3}},f_{{2,2}},f_{{2,3}},f_{{3,3}},f_{{2,2,2}},f_{{2,2,3}},f_{{2,3,3}},f_{{3,3,3}},f_{{2,2,2,2}},f_{{2,2,2,3}},f_{{2,2,3,3}},f_{{2,3,3,3}},f_{{3,3,3,3}},\\
f_{{2,2,2,2,2}},f_{{2,2,2,2,3}},f_{{2,2,2,3,3}},f_{{2,2,3,3,3}},f_{{2,3,3,3,3}},f_{{3,3,3,3,3}},d_{{1,1}},d_{{1,2}},d_{{1,3}},d_{{1,4}},d_{{1,5}},d_{{1,6}},d_{{1,7}},d_{{1,8}},d_{{1,9}},d_{{1,10}},\\
d_{{1,11}},d_{{1,12}},d_{{1,13}},d_{{1,14}},d_{{1,15}},d_{{2,1}},d_{{2,2}},d_{{2,3}},d_{{2,4}},d_{{2,5}},d_{{2,6}},d_{{2,7}},d_{{2,8}},d_{{2,9}},d_{{2,10}},d_{{3,1}},d_{{3,2}},d_{{3,3}},d_{{3,4}},\\
d_{{3,5}},d_{{3,6}},d_{{4,1}},d_{{4,2}},d_{{4,3}},d_{{5,1}},K, L, M)=0\\
\end{array}
\end{equation}
\end{footnotesize}
gives the implicit solution of the variables $K, L$ and $M$ with
the differential invariants.

In this paper, we are interested in finding the fifth order
differential invariants and differential operators of the general
class $y''=f(x,y,y')$  under point transformations (\ref{a9}). So,
according to the theory of invariants of infinite transformation
groups \cite{Ovsiannikov}, the invariant criterion (\ref{ccccc})
should be split by the functions  $\xi(x,y)$ and $\eta(x, y)$ and
their derivatives. This gives rise to a homogeneous linear system
of partial differential equations (PDEs):
\begin{equation}
\begin{array}{ll}\label{c3}
X_i J=0,~T_i J=0,~i=1\dots36,\\
\end{array}
\end{equation}
where $X_i,  i=1\dots36,$ are the differential operators
corresponding to the coefficients of the following derivatives of
$\eta(x, y)$ up to the seven order in the invariant criterion
\begin{footnotesize}
\begin{equation}
\begin{array}{ll}\label{c4}
\eta_{{}},\eta_{{1}},\eta_{{2}},\eta_{{1,1}},\eta_{{1,2}},\eta_{{2,2}},\eta_{{1,1,1}},\eta_{{1,1,2}},\eta_{{1,2,2}},\eta_{{2,2,2}},\eta_{{1,1,1,1}},\eta_{{1,1,1,2}},\eta_{{1,1,2,2}},\eta_{{1,2,2,2}},\eta_{{2,2,2,2}},\eta_{{1,1,1,1,1}},\eta_{{1,1,1,1,2}},\\
\eta_{{1,1,1,2,2}},\eta_{{1,1,2,2,2}},\eta_{{1,2,2,2,2}},\eta_{{2,2,2,2,2}},\eta_{{1,1,1,1,1,1}},\eta_{{1,1,1,1,1,2}},\eta_{{1,1,1,1,2,2}},\eta_{{1,1,1,2,2,2}},\eta_{{1,1,2,2,2,2}},\eta_{{1,2,2,2,2,2}},\eta_{{2,2,2,2,2,2}},\\
\eta_{{1,1,1,1,1,1,1}},\eta_{{1,1,1,1,1,1,2}},\eta_{{1,1,1,1,1,2,2}},\eta_{{1,1,1,1,2,2,2}},\eta_{{1,1,1,2,2,2,2}},\eta_{{1,1,2,2,2,2,2}},\eta_{{1,2,2,2,2,2,2}},\eta_{{2,2,2,2,2,2,2}}\\
\end{array}
\end{equation}
\end{footnotesize}
and $T_i, i=1\dots36,$ are the differential operators
corresponding to the coefficients of the following derivatives of
$\xi(x,y)$ up to the seven order  in the invariant criterion
\begin{footnotesize}
\begin{equation}
\begin{array}{ll}\label{c5}
\xi_{{}},\xi_{{1}},\xi_{{2}},\xi_{{1,1}},\xi_{{1,2}},\xi_{{2,2}},\xi_{{1,1,1}},\xi_{{1,1,2}},\xi_{{1,2,2}},\xi_{{2,2,2}},\xi_{{1,1,1,1}},\xi_{{1,1,1,2}},\xi_{{1,1,2,2}},\xi_{{1,2,2,2}},\xi_{{2,2,2,2}},\xi_{{1,1,1,1,1}},\xi_{{1,1,1,1,2}},\\
\xi_{{1,1,1,2,2}},\xi_{{1,1,2,2,2}},\xi_{{1,2,2,2,2}},\xi_{{2,2,2,2,2}},\xi_{{1,1,1,1,1,1}},\xi_{{1,1,1,1,1,2}},\xi_{{1,1,1,1,2,2}},\xi_{{1,1,1,2,2,2}},\xi_{{1,1,2,2,2,2}},\xi_{{1,2,2,2,2,2}},\xi_{{2,2,2,2,2,2}},\\
\xi_{{1,1,1,1,1,1,1}},\xi_{{1,1,1,1,1,1,2}},\xi_{{1,1,1,1,1,2,2}},\xi_{{1,1,1,1,2,2,2}},\xi_{{1,1,1,2,2,2,2}},\xi_{{1,1,2,2,2,2,2}},\xi_{{1,2,2,2,2,2,2}},\xi_{{2,2,2,2,2,2,2}}\\
\end{array}
\end{equation}
\end{footnotesize}
Functionally independent solutions of the system (\ref{c3})
provide all independent differential invariants of $y''=f(x,y,y')$
up to the fifth order under point transformations, as well as an
implicit solution of the variables $K, L$ and $M$ which provide
the differential operators via (\ref{cc}). The solution of system
(\ref{c3}) is found in many steps using Maple as follows:

First, let us consider the subsystem induced by the sixth and
seventh derivatives of $\xi$ and $\eta$
\begin{equation}\label{c10}
X_i J=0,~T_i J=0,~i=22\dots36.
\end{equation}
where the operators $X_i$ and $T_i,~i=22\dots36$ are given in
Appendix A in term of the variables $z_i, i=1\dots62$  after
relabeling the variables
\begin{small}
\begin{equation}\label{c6}
\begin{array}{ll}
x,y_{{}},y_{{1}},f_{{}},f_{{2}},f_{{3}},f_{{2,2}},f_{{2,3}},f_{{3,3}},f_{{2,2,2}},f_{{2,2,3}},f_{{2,3,3}},f_{{3,3,3}},f_{{2,2,2,2}},f_{{2,2,2,3}},f_{{2,2,3,3}},f_{{2,3,3,3}},f_{{3,3,3,3}},\\
f_{{2,2,2,2,2}},f_{{2,2,2,2,3}},f_{{2,2,2,3,3}},f_{{2,2,3,3,3}},f_{{2,3,3,3,3}},f_{{3,3,3,3,3}},d_{{1,1}},d_{{1,2}},d_{{1,3}},d_{{1,4}},d_{{1,5}},d_{{1,6}},d_{{1,7}},d_{{1,8}},d_{{1,9}},d_{{1,10}},\\
d_{{1,11}},d_{{1,12}},d_{{1,13}},d_{{1,14}},d_{{1,15}},d_{{2,1}},d_{{2,2}},d_{{2,3}},d_{{2,4}},d_{{2,5}},d_{{2,6}},d_{{2,7}},d_{{2,8}},d_{{2,9}},d_{{2,10}},d_{{3,1}},d_{{3,2}},d_{{3,3}},d_{{3,4}},\\
d_{{3,5}},d_{{3,6}},d_{{4,1}},d_{{4,2}},d_{{4,3}},d_{{5,1}},K, L, M\\
\end{array}
\end{equation}
\end{small}
by the variables $z_i, i=1\dots62$, respectively.

In 62-dimensional space of variables $z_i, i=1\dots62$, the rank
of the system (\ref{c10}) is 16, so it has 46 functionally
independent solutions which are given as:
\begin{footnotesize}
\begin{equation}\label{c11}
\begin{array}{ll}
l_{{1}}=z_{{1}},l_{{2}}=z_{{2}},l_{{3}}=z_{{3}},l_{{4}}=z_{{4}},l_{{5}}=z_{{5}},l_{{6}}=z_{{6}},l_{{7}}=z_{{7}},l_{{8}}=z_{{8}},l_{{9}}=z_{{9}},l_{{10}}=z_{{10}},l_{{11}}=z_{{11}},l_{{12}}=z_{{12}},l_{{13}}=z_{{13}},\\
l_{{14}}=z_{{15}},l_{{15}}=z_{{16}},l_{{16}}=z_{{17}},l_{{17}}=z_{{18}},l_{{18}}=z_{{21}},l_{{19}}=z_{{22}},l_{{20}}=z_{{23}},l_{{21}}=z_{{24}},l_{{22}}=z_{{25}},l_{{23}}=z_{{26}},l_{{24}}=z_{{27}},\\
l_{{25}}=z_{{28}},l_{{26}}=z_{{29}},l_{{27}}=z_{{30}},l_{{28}}=z_{{32}},l_{{29}}=z_{{33}},l_{{30}}=z_{{34}},l_{{31}}=z_{{37}},l_{{32}}=z_{{38}},l_{{33}}=z_{{39}},l_{{34}}=z_{{40}},l_{{35}}=z_{{41}},\\
l_{{36}}=z_{{42}},l_{{37}}=z_{{44}},l_{{38}}=z_{{45}},l_{{39}}=z_{{48}},l_{{40}}=z_{{49}},l_{{41}}=z_{{50}},l_{{42}}=z_{{52}},l_{{43}}=z_{{55}},l_{{44}}=z_{{60}},l_{{45}}=z_{{61}},l_{{46}}=z_{{62}}\\
\end{array}
\end{equation}
\end{footnotesize}
Second, let us consider the subsystem induced by the fifth
derivatives of $\xi$ and $\eta$
\begin{equation}\label{c13}
X_i J=0,~ T_i J=0,~i=16\dots21.
\end{equation}
where the inherited operators $X_i$ and $T_i,~i=16\dots21$ are
given in Appendix B in term of the new variables $l_i,
i=1\dots46$.

In 46-dimensional space of variables $l_i, i=1\dots46$, the rank
of the system (\ref{c13}) is 10, so it has 36 functionally
independent solutions which are given as:
\begin{equation}\label{c14}
\begin{array}{llllll}
m_{{1}}=l_{{1}},m_{{2}}=l_{{2}},m_{{3}}=l_{{3}},m_{{4}}=l_{{4}},m_{{5}}=l_{{5}},m_{{6}}=l_{{6}},m_{{7}}=l_{{7}},m_{{8}}=l_{{8}},m_{{9}}=l_{{9}},\\
m_{{10}}=l_{{11}},m_{{11}}=l_{{12}},m_{{12}}=l_{{13}},m_{{13}}=l_{{15}},m_{{14}}=l_{{16}},m_{{15}}=l_{{17}},m_{{16}}=l_{{19}},m_{{17}}=l_{{20}},\\
m_{{18}}=l_{{21}},m_{{19}}=l_{{22}},m_{{20}}=l_{{23}},m_{{21}}=l_{{24}},m_{{22}}=l_{{26}},m_{{23}}=l_{{27}},m_{{24}}=l_{{29}},m_{{25}}=l_{{30}},\\
m_{{26}}=l_{{32}},m_{{27}}=l_{{33}},m_{{28}}=l_{{34}},m_{{29}}=l_{{36}},m_{{30}}=l_{{38}},m_{{31}}=-4\,l_{{28}}+6\,l_{{10}}+l_{{39}},m_{{32}}=l_{{40}},\\
m_{{33}}=-4\,l_{{37}}+6\,l_{{25}}+l_{{43}},m_{{34}}=l_{{44}},m_{{35}}=l_{{45}},m_{{36}}=l_{{46}}
\end{array}
\end{equation}
Third, let us consider the subsystem induced by the fourth
derivatives of $\xi$ and $\eta$
\begin{equation}\label{c15}
X_i J=0,~ T_i J=0,~i=11\dots15.
\end{equation}
where the inherited operators $X_i$ and $T_i,~i=11\dots15$ are
given in Appendix C in term of the new variables $m_i,
i=1\dots36$.

In 36-dimensional space of variables $m_i, i=1\dots36$, the rank
of the system (\ref{c15}) is 10, so it has 26 functionally
independent solutions which are given as:
\begin{equation}\label{c16}
\begin{array}{llllll}
n_{{1}}=m_{{1}},n_{{2}}=m_{{2}},n_{{3}}=m_{{3}},n_{{4}}=m_{{4}},n_{{5}}=m_{{5}},n_{{6}}=m_{{6}},n_{{7}}=m_{{8}},n_{{8}}=m_{{9}},n_{{9}}=m_{{11}},\\
n_{{10}}=m_{{12}},n_{{11}}=m_{{14}},n_{{12}}=m_{{15}},n_{{13}}=m_{{17}},n_{{14}}=m_{{18}},n_{{15}}=m_{{19}},n_{{16}}=m_{{21}},n_{{17}}=m_{{23}},\\
n_{{18}}=m_{{25}},n_{{19}}=m_{{27}},n_{{20}}=6\,m_{{7}}-4\,m_{{22}}+m_{{30}},n_{{21}}=-3\,m_{{9}}m_{{7}}+m_{{12}}m_{{20}}-m_{{6}}m_{{24}}+\\
4\,m_{{6}}m_{{10}}+m_{{31}},n_{{22}}=-2\,m_{{24}}+2\,m_{{10}}+m_{{32}},n_{{23}}=6\,m_{{6}}m_{{7}}-3\,m_{{9}}m_{{20}}+m_{{33}},\\
n_{{24}}=m_{{34}},n_{{25}}=m_{{35}},n_{{26}}=m_{{36}}\\
\end{array}
\end{equation}
Fourth, let us consider the subsystem induced by the third
derivatives of $\xi$ and $\eta$
\begin{equation}\label{c17}
X_i J=0,~ T_i J=0,~i=7\dots10.
\end{equation}
where the inherited operators $X_i$ and $T_i,~i=7\dots10$ are
given in Appendix D in term of the new variables $n_i,
i=1\dots26$.

In 26-dimensional space of variables $n_i, i=1\dots26$, the rank
of the system (\ref{c17}) is 8, so it has 18 functionally
independent solutions which are given as:
\begin{equation}\label{c18}
\begin{array}{llllll}
t_{{1}}=n_{{1}},t_{{2}}=n_{{2}},t_{{3}}=n_{{3}},t_{{4}}=n_{{4}},t_{{5}}=n_{{6}},t_{{6}}=n_{{8}},t_{{7}}=n_{{10}},t_{{8}}=n_{{12}},t_{{9}}=n_{{13}},t_{{10}}=n_{{14}},\\
t_{{11}}=n_{{19}},t_{{12}}=-3\,n_{{8}}n_{{5}}-n_{{6}}n_{{17}}+4\,n_{{7}}n_{{6}}+n_{{20}},t_{{13}}=4\,{n_{{7}}}^{2}-n_{{17}}n_{{7}}+2\,n_{{5}}n_{{18}}-6\,n_{{5}}n_{{9}}+n_{{21}},\\
t_{{14}}=-2\,n_{{10}}n_{{5}}-n_{{8}}n_{{17}}+n_{{7}}n_{{8}}+n_{{16}}n_{{10}}+n_{{6}}n_{{18}}+n_{{22}},t_{{15}}=-3\,n_{{6}}n_{{8}}n_{{5}}-n_{{17}}{n_{{6}}}^{2}-3\,n_{{17}}n_{{5}}\\
+4\,{n_{{6}}}^{2}n_{{7}}+4\,n_{{16}}n_{{7}}-n_{{16}}n_{{17}}+n_{{23}},t_{{16}}=n_{{24}},t_{{17}}=n_{{25}},t_{{18}}=n_{{26}}\\
\end{array}
\end{equation}
Finally, let us consider the subsystem induced by the zero, first
and second derivatives of $\xi$ and $\eta$
\begin{equation}\label{c19}
X_i J=0,~ T_i J=0,~i=1\dots6.
\end{equation}
where the inherited operators $X_i$ and $T_i,~i=1\dots6$ are given
in Appendix E in term of the new variables $t_i, i=1\dots18$ which
can be rewritten, by backing substitution, as
\begin{equation}
\begin{array}{llllll}
t_{{1}}=x,t_{{2}}=y_{{}},t_{{3}}=y_{{1}},t_{{4}}=f_{{}},t_{{5}}=f_{{3}},t_{{6}}=f_{{3,3}},t_{{7}}=f_{{3,3,3}},t_{{8}}=f_{{3,3,3,3}},\\
t_{{9}}=f_{{2,3,3,3,3}},t_{{10}}=f_{{3,3,3,3,3}},t_{{11}}=\dot{D}_x{f_{3,3,3,3}},t_{{16}}=K_{{}},t_{{17}}=L_{{}},t_{{18}}=M_{{}}\\
\end{array}
\end{equation}
and
\begin{equation}
\begin{array}{llllll}
t_{{12}}=&4\,f_{{3}}f_{{2,3}}-f_{{3}}\dot{D}_x{f_{3,3}}-3\,f_{{3,3}}f_{{2}}+6\,f_{{2,2}}+\dot{D}_x^2{f_{3,3}}-4\,\dot{D}_x{f_{2,3}},\\
t_{{13}}=&\tilde{D}_y  {t_{12}}\\
t_{{14}}=&\tilde{D}_p  {t_{12}},\\
t_{{15}}=&f_3~{t_{12}} +\dot{D}_x {t_{12}}\\
\end{array}
\end{equation}
It should be noted here that $t_8$ and $t_{12}$ are the fourth
order Tress´e relative invariants. It is well known that a scalar
second-order ODE is linearizable via a point transformation if and
only if they both vanish identically as shown by Theorem 1.1.
Moreover, it is noted that $t_{{13}}, t_{{14}}$ and $t_{{15}}$
vanish identically when $t_{{12}}=0$.

One can see  that the operators $X_i$ and $T_i,~i=1\dots6$ form a
Lie algebra $\mathcal{L}_{12}$ with the nonzero commutators
\begin{equation}
\begin{tabular}{llllll}
$[X_2, X_3] = X_2$,& $[X_2, X_5] = 2~X_4$,& $[X_2, X_6] = X_5$,& $[X_2, T_2] =-X_2$,\\
$[X_2, T_3]=T_2-X_3$,&$[X_2, T_4]=-X_4$,&$[X_2, T_5] =-X_5+2~T_4$,& $[X_2, T_6] = -X_6+T_5$,\\
$[X_3, X_4]=-X_4$,&$[X_3, X_6] =X_6$,&$[X_3, T_3] = T_3$,&$[X_3, T_5] = T_5$,\\
$[X_3, T_6] = 2~T_6$,&$[X_4, T_2]=-2~X_4$,& $[X_4, T_3] = T_4-X_5$,&$[X_5, T_2] = -X_5$,\\
$[X_5,T_3] = -2~X_6+T_5$,&$[X_6, T_3] = T_6$,& $[T_2, T_3]=-T_3$,& $[T_2, T_4] =T_4$,\\
$[T_2, T_6] = -T_6$,&$[T_3, T_4] =T_5$,& $[T_3, T_5] = 2~T_6$\\
\end{tabular}
\end{equation}
Moreover, the projection of the operators $X_i$ and
$T_i,~i=1\dots6$ on the 4-dimensional space of variables $t_i,
i=1\dots4$ are the generators of the original infinite Lie algebra
spanned by the infinitesimal operators (\ref{a11}) before the
prolongation to the fifth order.

In section 4, the joint invariants  of the operators (\ref{c19})
provide all differential invariants of $y''=f(x,y,y')$, with
$f_{3,3,3,3}\ne0$, up to the fifth order under point
transformations.
\section{The infinitesimal point equivalence transformations} In
order to find continuous group of equivalence transformations of
the class (\ref{w1}) we consider the arbitrary function $f$ that
appears in our equation as a dependent variable and the variables
$x,y,y'=p$ as independent variables and apply the Lie
infinitesimal invariance criterion \cite{Ovsiannikov}, that is we
look for the infinitesimal $\xi, \eta$ and $\mu$ of the
equivalence operator $Y$:
\begin{equation}
Y=\xi(x,y)\partial_x+\eta(x,y)\partial_y+\mu(x,y,p,f)\partial_f,\\
\end{equation}
such that its prolongation leaves the equation (\ref{w1})
invariant.

The prolongation of operator $Y$ can be given using (\ref{a5}) as
\begin{equation}
Y=\xi(x,y)D_x+W\partial_y+D_x(W)\partial_{p}+D_x^2(W)\partial_{y''}+\mu(x,y,p,f)\partial_f,\\
\end{equation}
where $$ D_x  = \frac{\partial } {{\partial x}} + p\frac{\partial
}{{\partial y}} + y''\frac{\partial } {{\partial p}}
+y'''\frac{\partial} {{\partial y''}}+...$$ is the operator of
total derivative and $W = \eta(x,y) - \xi(x,y) p$ is the
characteristic of infinitesimal operator
$X=\xi(x,y)\partial_x+\eta(x,y)\partial_y$.

So, the Lie infinitesimal invariance criterion gives
$\mu=\dot{D}_x^2(W)+\xi(x,y)\dot{D}_x f$ for arbitrary functions
$\xi(x,y)$ and $\eta(x,y)$ where $ \dot{D}_x = \frac{\partial }
{{\partial x}} + p\frac{\partial }{{\partial y}} +f\frac{\partial}
{{\partial p}}$.

Thus, equation (\ref{w1}) admits an infinite continuous group of
equivalence transformations generated by the Lie algebra
$\mathcal{L}_\mathcal{E}$ spanned by the following infinitesimal
operators
\begin{equation}\label{e3}
U=\xi(x,y)\frac{\partial } {{\partial x}}-p D_x(\xi)\partial_{p}-(2 f D_x(\xi)+p \dot{D}_x^2(\xi))\partial_f,\\
\end{equation}
\begin{equation}\label{e4}
V=\eta(x,y)\partial_y+D_x(\eta)\partial_{p}+\dot{D}_x^2(\eta)\partial_f,\\
\end{equation}
The infinitesimal point equivalence transformations
(\ref{e3})-(\ref{e4}) can be written in the finite form as in
(\ref{a12})-(\ref{a13}), respectively, where $\phi$ and $\psi$ are
arbitrary functions of the indicated variables.
\section{Fifth-order differential invariants and invariant
equations under point transformations} In this section, we derive
all the fifth-order differential invariants of the general class
\\$y''=f(x,y,y')$, with $f_{3,3,3,3}\ne0$,  under point
transformations (\ref{a9}). Moreover, the invariant
differentiation operators \cite{Ovsiannikov} are constructed in
order to get some higher-order differential invariants from the
lower-order ones. Precisely, we obtain the following theorem.
\begin{theorem}
Every second-order ODE $y''=f(x,y,y')$, with $f_{3,3,3,3}\ne0$,
belongs to one of two classes of equations. For the first class of
equation ($\nu_1 \ne 0$), there are three fifth order differential
invariants, under point transformations,
\begin{equation}\label{d1}
\begin{array}{lllll}
\beta_{{1}}={{\nu_{{2}}}^{4}}{{\nu_{{1}}}^{-\frac{7}{2}}},&\beta_{{2}}= {{\nu_{{3}}}^{4}}{{\nu_{{1}}}^{-\frac{11}{2}}},&\beta_{{3}}={{\nu_{{4}}}^{4}}{{\nu_{{1}}}^{-5}},\\
\end{array}
\end{equation}
and three invariant differential operators
\begin{footnotesize}
\begin{equation}\label{d2}
\begin{array}{llll}
\mathcal{D}_1={(f_{{3,3,3,3}})}^{-\frac{2}{5}}\nu_1^{\frac{1}{8}}\tilde{D}_p,\\
\mathcal{D}_2={(f_{{3,3,3,3}})}^{\frac{1}{5}}\nu_1^{-\frac{3}{8}}\left(\tilde{D}_x+p~\tilde{D}_y+f\tilde{D}_p\right),\\
\mathcal{D}_3={(f_{{3,3,3,3}})}^{-\frac{6}{5}}\nu_1^{-\frac{1}{4}}\left(f_{{3,3,3,3,3}}\tilde{D}_x+(5f_{3, 3, 3, 3}+p~f_{3, 3, 3, 3, 3})\tilde{D}_y+(10f_3 f_{3, 3, 3, 3}+ff_{3, 3, 3, 3, 3}+5~\dot{D}_x{f_{3,3,3,3}})\tilde{D}_p\right),\\
\end{array}
\end{equation}
\end{footnotesize}
which satisfy the higher order relations
\begin{equation}\label{ast1}
\begin{array}{llll}
\mathcal{D}_1\mathcal{D}_2H-\mathcal{D}_2\mathcal{D}_1 H-\rho_1\mathcal{D}_1H-\rho_2\mathcal{D}_2 H-\rho_3\mathcal{D}_3 H=0,\\
\mathcal{D}_1\mathcal{D}_3H-\mathcal{D}_3\mathcal{D}_1H-\sigma_1\mathcal{D}_1H-\sigma_2\mathcal{D}_2H-\sigma_3\mathcal{D}_3H=0,\\
\mathcal{D}_2\mathcal{D}_3H-\mathcal{D}_3\mathcal{D}_2H-\omega_1\mathcal{D}_1H-\omega_2\mathcal{D}_2H-\omega_3\mathcal{D}_3H=0,\\
\end{array}
\end{equation}
for any differential invariant $H$.

However, there is no fifth-order differential invariants for the
second class ($\nu_1=0$), where $\nu_1, \nu_2, \nu_3$ and $\nu_4$
are the relative invariants given by (\ref{c20}) and the
commutator invariants $\rho_1, \rho_2, \rho_3$, $\sigma_1,
\sigma_2, \sigma_3$ and $\omega_1, \omega_2, \omega_3$ can be
given by (\ref{ast3}).
\end{theorem}
\proof The joint invariants  of the operators (\ref{c19}) provide
all differential invariants of $y''=f(x,y,y')$, with
$f_{3,3,3,3}\ne0$, up to the fifth order under point
transformations, as well as an implicit solution of the variables
$K, L$ and $M$ which provide the differential operators via
(\ref{cc}).

The joint invariants of the first derived subgroup of
$\mathcal{L}_{12}$ can be given for the case $f_{3,3,3,3}\ne0$
after backing substitution as an arbitrary function
$J(x,y,\nu_{{1}},\nu_{{2}},\nu_{{3}},\nu_{{4}}.\nu_{{5}},\nu_{{6}},\nu_{{7}}),$
where
\begin{footnotesize}
\begin{equation}\label{c20}
\begin{array}{lllll}
\nu_{{1}}=&{t_{{8}}}^{\frac{1}{5}}t_{{12}}, \\
\nu_{{2}}=&{t_{{8}}}^{-\frac{6}{5}}\left(t_{{12}}t_{{10}}+5\,t_{{8}}t_{{14}}\right),\\
\nu_{{3}}=&{t_{{8}}}^{-\frac{3}{5}}\left(  5\,t_{{11}}t_{{12}}+\left( 7\,t_{{5}}t_{{12}}+t_{{15}} \right) t_{{8} }\right),\\
\nu_{{4}}=&{t_{{8}}}^{-2}\left(\left(10\,t_{{5}}t_{{14}}-5\,t_{{12}}t_{{6}}+5\,t_{{13}} \right) {t_{{8}}}^{2}+ \left(  \left( -3\,t_{{5}}t_{{10}}-5\,t_{{9}} \right)t_{{ 12}}+t_{{15}}t_{{10}}+5\,t_{{14}}t_{{11}} \right) t_{{8}}\right)\\
\end{array}
\end{equation}
\end{footnotesize}
and
\begin{equation}\label{c21}
\begin{array}{llllll}
\nu_5={f_{{3,3,3,3}}}^{\frac{1}{5}} \left(L_{{}}- K_{{}}y_{{1}} \right), \\
\nu_6=\frac{1}{5}{f_{{3,3,3,3}}}^{-\frac{6}{5}} \left(5\,K_{{}}f_{{3,3,3,3}}+K_{{}}y_{{1}}f_{{3,3,3,3,3}}-f_{{3,3,3,3,3}}L_{{}}\right),\\
\nu_7={f_{{3,3,3,3}}}^{-\frac{3}{5}}\left(2\,f_{{3,3,3,3}}K_{{}}f_{{3}}y_{{1}}-f_{{3,3,3,3}}K_{{}}f_{{}}-2\,f_{{3,3,3,3}}f_{{3}}L_{{}}+f_{{3,3,3,3}}M_{{}}+\dot{D}_x{f_{3,3,3,3}}(K_{{}}y_{{1}}-L_{{}})\right).\\
\end{array}
\end{equation}
The non-zero inheritance of the operators $X_i$ and
$T_i,~i=1\dots12$ in term of the new variables $x,y,\nu_i,
i=1\dots7$ is
\begin{equation}\label{c22}
\begin{array}{llllll}
X_1=&[0,1,0,0,0,0,0,0,0],\\
T_1=&[1,0,0,0,0,0,0,0,0],\\
X_3=&[0,0,-\frac{8}{5}\,\nu_{{1}},-\frac{7}{5}\,\nu_{{2}},-\frac{11}{5}\,\nu_{{3}},-2\,\nu_{{4}},\frac{2}{5}\,\nu_{{5}},\frac{3}{5}\,\nu_{{6}},-\frac{1}{5}\nu_{{7}}],\\
T_2=&X_3.\\
\end{array}
\end{equation}
The joint invariants of the operators (\ref{c22}) are the
invariants of the operator
\begin{equation}\label{c23}
\begin{array}{llllll}
Z=&8\,\nu_{{1}}\frac{\partial}{\partial \nu_1}+7\,\nu_{{2}}\frac{\partial}{\partial \nu_2}+11\,\nu_{{3}}\frac{\partial}{\partial \nu_3}+10\,\nu_{{4}}\frac{\partial}{\partial \nu_4}-2\,\nu_{{5}}\frac{\partial}{\partial \nu_5}-3\,\nu_{{6}}\frac{\partial}{\partial \nu_6}+\nu_{{7}}\frac{\partial}{\partial \nu_7}.\\
\end{array}
\end{equation}
The invariants of the operators (\ref{c23}) can be given using
characteristic method for two classes as follows:

(1) First class of equation ($\nu_1 \ne 0$)
\begin{equation}\label{cc20}
\begin{array}{llllll}
\beta_{{1}}={{\nu_{{2}}}^{4}}{{\nu_{{1}}}^{-\frac{7}{2}}},&\beta_{{2}}= {{\nu_{{3}}}^{4}}{{\nu_{{1}}}^{-\frac{11}{2}}},&\beta_{{3}}={{\nu_{{4}}}^{4}}{{\nu_{{1}}}^{-5}},\\
\end{array}
\end{equation}
and
\begin{equation}\label{cc21}
\begin{array}{llllll}
\gamma_1={\nu_{{5}}}^{8}{\nu_{{1}}}^{2},~\gamma_2={\nu_{{6}}}^{8}{\nu_{{1}}}^{3},~\gamma_3={\frac {{\nu_{{7}}}^{8}}{\nu_{{1}}}}\\
\end{array}
\end{equation}
(2) Second class of equation ($\nu_1=0$) does not have fifth-order
differential invariants independent from the variables $K, L$ and
$M$. This because of vanishing the variables $t_{{13}}, t_{{14}}$
and $t_{{15}}$ identically when $t_{{12}}=0$, and so
$\nu_2=\nu_3=\nu_4=0$ whenever $\nu_1=0$.

Regarding  the invariant differentiation operators, $\gamma_1,
\gamma_2$ and $\gamma_3$ are the only invariants depending on the
variables $K, L$ and $M$. Then the general solution of (\ref{ccc})
can be given implicitly as
\begin{equation}\label{cc26}
\begin{array}{llll}
\gamma_1={\it F_1},~\gamma_2={\it F_2},~\gamma_3={\it F_3} ,\\
\end{array}
\end{equation}
where $F_1, F_2$ and $F_3$ are the arbitrary functions of
differential invariants $\beta_{{i}},~i=1\dots3.$

Solving system (\ref{cc26}) gives the variables $K, L$ and $M$ in
terms of three arbitrary functions $F_1, F_2$ and $F_3$ which
provide three independent invariant differentiation operators
$\mathcal{D}_1,\mathcal{D}_2$ and $\mathcal{D}_3$ via (\ref{cc}).

Finally, since the matrix \begin{equation}
\begin{array}{cc}
A=\left(%
\begin{array}{ccc}
\mathcal{D}_1 x&\mathcal{D}_2 x&\mathcal{D}_3 x\\
\mathcal{D}_1 y&\mathcal{D}_2 y&\mathcal{D}_3 y\\
\mathcal{D}_1 p&\mathcal{D}_2 p&\mathcal{D}_3 p\\
\end{array}%
\right)
\end{array}
\end{equation}
is an invertible matrix with the non-zero determinant
$J=5~{f_{{3,3,3,3}}}^{-\frac{2}{5}}{\nu_1}^{-\frac{1}{2}}$, then
the invariant differential operators should satisfy the
commutation relations
\begin{equation}\label{ast2}
\begin{tabular}{llll}
$[\mathcal{D}_1,\mathcal{D}_2]=\rho_1\mathcal{D}_1+\rho_2\mathcal{D}_2+\rho_3\mathcal{D}_3$,\\
$[\mathcal{D}_1,\mathcal{D}_3]=\sigma_1\mathcal{D}_1+\sigma_2\mathcal{D}_2+\sigma_3\mathcal{D}_3$,\\
$[\mathcal{D}_2,\mathcal{D}_3]=\omega_1\mathcal{D}_1+\omega_2\mathcal{D}_2+\omega_3\mathcal{D}_3$,\\
\end{tabular}
\end{equation}
where
\begin{equation}\label{ast3}
\begin{array}{cc}
\left(%
\begin{array}{cccc}
\rho_1\\
\rho_2\\
\rho_3\\
\end{array}%
\right) =
A^{-1}
\left(%
\begin{array}{cccc}
\mathcal{D}_1\mathcal{D}_2x-\mathcal{D}_2\mathcal{D}_1x\\
\mathcal{D}_1\mathcal{D}_2y-\mathcal{D}_2\mathcal{D}_1y\\
\mathcal{D}_1\mathcal{D}_2p-\mathcal{D}_2\mathcal{D}_1p\\
\end{array}%
\right), &
\left(%
\begin{array}{cccc}
\sigma_1\\
\sigma_2\\
\sigma_3\\
\end{array}%
\right) =A^{-1}
\left(%
\begin{array}{cccc}
\mathcal{D}_1\mathcal{D}_3x-\mathcal{D}_3\mathcal{D}_1x\\
\mathcal{D}_1\mathcal{D}_3y-\mathcal{D}_3\mathcal{D}_1y\\
\mathcal{D}_1\mathcal{D}_3p-\mathcal{D}_3\mathcal{D}_1p\\
\end{array}%
\right),\\
\left(%
\begin{array}{cccc}
\omega_1\\
\omega_2\\
\omega_3\\
\end{array}%
\right) =A^{-1}
\left(%
\begin{array}{cccc}
\mathcal{D}_2\mathcal{D}_3x-\mathcal{D}_3\mathcal{D}_2x\\
\mathcal{D}_2\mathcal{D}_3y-\mathcal{D}_3\mathcal{D}_2y\\
\mathcal{D}_2\mathcal{D}_3p-\mathcal{D}_3\mathcal{D}_2p\\
\end{array}%
\right).
\end{array}
\end{equation}
Hence the commutator identities (\ref{ast2}) can be applied to any
differential invariants $H$ to give the higher order relations
(\ref{ast1}).
\endproof
\section{Application}
In this section, invariant description of all the canonical forms
in the complex plane  \cite{Olver1995} for second-order ODEs
$y''=f(x,y,y')$ where both of the two Tress\'e relative invariants
(\ref{wer}) are non-zero is provided. Moreover, one example of the
second class ($\nu_1=0$) is given from the canonical forms of
second order ODE in the real plane \cite{Mahomed2007}.
\begin{example}\rm \label{ex2}
 Consider the canonical form of second order ODE in the complex
plane  with three infinitesimal symmetries \cite{Olver1995}
\begin{equation}
y''=c~\exp(-y').
\end{equation}
It is an equation of the first class ($\nu_1\ne0$), with the three
constant fifth-order differential invariants
\begin{equation}
\begin{array}{llllll}
\beta_1=-65536,~\beta_2=-65536,~\beta_3=2825761.\\
\end{array}
\end{equation}
\end{example}
\begin{example}\rm \label{ex3}
 Consider the canonical form of second order ODE in the complex
plane  with three infinitesimal symmetries \cite{Olver1995}
\begin{equation}
y''=c~y'^{(\frac{\alpha-2}{\alpha-1})},~ \alpha \ne 0,\frac{1}{2},1,2.\\
\end{equation}
It is an equation of the first class ($\nu_1\ne0$), with the three
fifth-order differential invariants
\begin{equation}
\begin{array}{llllll}
\beta_1=-4096\,{\frac { \left( \alpha+1\right)^{4}}{2\,{\alpha}^{3}-5\,{\alpha}^{2}+2\,\alpha}},~\beta_2=-4096\,{\frac { \left( \alpha+1\right) ^{4}}{2\,{\alpha}^{3}-5\,{ \alpha}^{2}+2\,\alpha}},~\beta_3={\frac { \left( 14\,{\alpha}^{2}+13\,\alpha+14 \right)^{4}}{{\alpha}^ {2} \left( 2\,\alpha-1 \right) ^{2} \left(\alpha-2 \right) ^{2}}}.\\
\end{array}
\end{equation}
As a special case, when  $\alpha=-1$, one have the second order
ODE
\begin{equation}
y''=c~y'^{\frac{3}{2}}.\\
\end{equation}
with the three fifth-order differential invariants
\begin{equation}
\begin{array}{llllll}
\beta_1=0,~\beta_2=0,~\beta_3=625.\\
\end{array}
\end{equation}
\end{example}
\begin{example}\rm \label{ex4}
 Consider the canonical form of second order ODE in the complex
plane  with three infinitesimal symmetries \cite{Olver1995}
\begin{equation}
y''=6~yy'-4~y^3+c~(y'-y^2)^{\frac{3}{2}},~c\ne\pm4i.
\end{equation}
It is an equation of the first class ($\nu_1\ne0$), with the three
fifth-order differential invariants
\begin{equation}
\begin{array}{llllll}
\beta_1=0,~\beta_2=0,~\beta_3=625\,{\frac {{c}^{2}}{16+{c}^{2}}}.\\
\end{array}
\end{equation}
\end{example}
\begin{example}\rm \label{ex5}
 Consider the canonical form of second order ODE in the complex
plane  with three infinitesimal symmetries \cite{Olver1995}
\begin{equation}
y''=6~yy'-4~y^3+c~(y'-y^2)^{\frac{3}{2}},~c=\pm4i.
\end{equation}
It is an equation of the second class ($\nu_1=0$), so it does not
have fifth-order differential invariants.
\end{example}
\begin{example}\rm \label{ex6}
 Consider the canonical form of second order ODE in the real
plane \cite{Mahomed2007}
\begin{equation}
x~y''= y'+{y'}^3+(1+y'^2)^{\frac{3}{2}}.
\end{equation}
It is an equation of the second class ($\nu_1=0$), so it does not
have fifth-order differential invariants.
\end{example}
\begin{example}\rm \label{ex7}
 Consider the canonical form of second order ODE in the complex
plane  with two infinitesimal symmetries \cite{Olver1995}
\begin{equation}
y''=f(y').
\end{equation}
It is an equation of the first class ($\nu_1\ne0$). It has three
non-constant fifth-order differential invariants. However, this
class can be characterized by the relation $\beta_1+\beta_2=0$ and
the Jacobian matrix
$\frac{\partial(\beta_1,\beta_2,\beta_3)}{\partial(x,y,p)}$ where
it has rank one.
\end{example}
\begin{example}\rm \label{ex8}
 Consider the canonical form of second order ODE in the complex
plane  with two infinitesimal symmetries \cite{Olver1995}
\begin{equation}
y''=y'+f(y'-y).
\end{equation}
For the case ($\nu_1\ne0$), it has three non-constant fifth-order
differential invariants. However, this class can be characterized
by the relation $\beta_1+\beta_2\ne0$ and the Jacobian matrix
$\frac{\partial(\beta_1,\beta_2,\beta_3)}{\partial(x,y,p)}$ where
it has rank one.
\end{example}
\begin{example}\rm \label{ex9}
 Consider the canonical form of second order ODE in the complex
plane  with one infinitesimal symmetries \cite{Olver1995}
\begin{equation}
y''=f(x,y').
\end{equation}
For the case ($\nu_1\ne0$), it has three non-constant fifth-order
differential invariants. However, this class can be characterized
by the Jacobian matrix
$\frac{\partial(\beta_1,\beta_2,\beta_3)}{\partial(x,y,p)}$ where
it has rank two.
\end{example}
\section{Conclusion}
The paper provides an extension of the work of Bagderina
\cite{Bagderina2013} who solved the equivalence problem for scalar
second-order ordinary differential equations (ODEs), cubic in the
first-order derivative, via point transformations. However, the
question is open for the general class $y''=f(x,y,y')$ which is
not cubic in the first-order derivative. Lie's infinitesimal
method was utilized to study the differential invariants of this
general class under an arbitrary point equivalence
transformations. All fifth order differential invariants and the
invariant differentiation operators were determined. These are
stated as Theorems 4.1 in Section 4.

As an application, the symmetry algebra of the second order ODE
$y''=f(x,y,y')$ where both of the two Tress\'e relative invariants
(\ref{wer}) are non-zero is characterized as follows:

1) The symmetry algebra is 3-dimensional iff the rank of the
Jacobian matrix
$\frac{\partial(\beta_1,\beta_2,\beta_3)}{\partial(x,y,p)}$ is
zero (the differential invariants $\beta_1, \beta_2$ and $\beta_3$
are constant).

2) The symmetry algebra is 2-dimensional iff the rank of the
Jacobian matrix
$\frac{\partial(\beta_1,\beta_2,\beta_3)}{\partial(x,y,p)}$ is
one.

3) The symmetry algebra is 1-dimensional iff the rank of the
Jacobian matrix
$\frac{\partial(\beta_1,\beta_2,\beta_3)}{\partial(x,y,p)}$ is
two.

Moreover, invariant description of all the canonical forms in the
complex plane for second-order ODEs $y''=f(x,y,y')$ where both of
the two Tress\'e relative invariants (\ref{wer}) are non-zero is
provided.
\subsection*{Acknowledgments}
 The author would like to thank the King Fahd University of
Petroleum and Minerals for its support and excellent research
facilities and I also want to thank Drs. Fazal Mahomed, Hassan
Azad and Tahir Mustafa for several discussions.
\section*{Appendix A: The differential operators of the homogeneous linear system of
PDEs (\ref{c10}) }
\begin{footnotesize}\tiny
\begin{eqnarray*}
X_{22} &=&[0,0,0,0,0,0,0,0,0,0,0,0,0,0,0,0,0,0,0,0,0,0,0,0,0,0,0,0,0,0,0,0,0,0,0,0,0,0,0,0,0,\\
&&0,0,0,0,0,0,0,0,0,0,0,0,0,0,1,0,0,0,0,0,0]\\
X_{23}&=&[0,0,0,0,0,0,0,0,0,0,0,0,0,0,0,0,0,0,0,0,0,0,0,0,0,0,0,0,0,0,0,0,0,0,0,0,0,0,0,0,0,0,\\
&&0,0,0,0,0,0,0,0,1,0,0,0,0,6\,z_{{3}},-z_{{6}},2,21\,z_{{4}},0,0,0]\\
X_{24}&=&[0,0,0,0,0,0,0,0,0,0,0,0,0,0,0,0,0,0,0,0,0,0,0,0,0,0,0,0,0,0,0,0,0,0,0,0,0,0,0,0,0,0,1,\\
&&0,0,0,0,0,0,0,5\,z_{{3}},0,-z_{{6}},2,0,15\,{z_{{3}}}^{2},15\,z_{{4}}-5\,z_{{6}}z_{{3}},10\,z_{{3}},105\,z_{{4}}z_{{3}},0,0,0]\\
X_{25}&=&[0,0,0,0,0,0,0,0,0,0,0,0,0,0,0,0,0,0,0,0,0,0,0,0,0,0,0,0,0,0,1,0,0,0,0,0,0,0,0,0,0,0,4\,z_{{3}},\\
&&0,0,-z_{{6}},2,0,0,0,10\,{z_{{3}}}^{2},0,10\,z_{{4}}-4\,z_{{6}}z_{{3}},8\,z_{{3}},0,20\,{z_{{3}}}^{3},60\,z_{{4}}z_{{3}}-10\,z_{{6}}{z_{{3}}}^{2},20\,{z_{{3}}}^{2},210\,z_{{4}}{z_{{3}}}^{2},0,0,0]\\
X_{26}&=&[0,0,0,0,0,0,0,0,0,0,0,0,0,1,0,0,0,0,0,0,0,0,0,0,0,0,0,0,0,0,3\,z_{{3}},0,0,0,-z_{{6}},2,0,0,0,0,0,0,6\,{z_{{3}}}^{2},\\
&&0,0,-3\,z_{{6}}z_{{3}}+6\,z_{{4}},6\,z_{{3}},0,0,0,10\,{z_{{3}}}^{3},0,30\,z_{{4}}z_{{3}}-6\,z_{{6}}{z_{{3}}}^{2},12\,{z_{{3}}}^{2},0,15\,{z_{{3}}}^{4},-10\,z_{{6}}{z_{{3}}}^{3}+90\,z_{{4}}{z_{{3}}}^{2},20\,{z_{{3}}}^{3},\\
&&210\,z_{{4}}{z_{{3}}}^{3},0,0,0]\\
X_{27}&=&[0,0,0,0,0,0,0,0,0,0,0,0,0,2\,z_{{3}},0,0,0,0,-z_{{6}},2,0,0,0,0,0,0,0,0,0,0,3\,{z_{{3}}}^{2},0,0,0,3\,z_{{4}}-2\,z_{{6}}z_{{3}},4\,z_{{3}},\\
&&0,0,0,0,0,0,4\,{z_{{3}}}^{3},0,0,12\,z_{{4}}z_{{3}}-3\,z_{{6}}{z_{{3}}}^{2},6\,{z_{{3}}}^{2},0,0,0,5\,{z_{{3}}}^{4},0,-4\,z_{{6}}{z_{{3}}}^{3}+30\,z_{{4}}{z_{{3}}}^{2},8\,{z_{{3}}}^{3},0,6\,{z_{{3}}}^{5},\\
&&60\,z_{{4}}{z_{{3}}}^{3}-5\,z_{{6}}{z_{{3}}}^{4},10\,{z_{{3}}}^{4},105\,z_{{4}}{z_{{3}}}^{4},0,0,0]\\
X_{28}&=&[0,0,0,0,0,0,0,0,0,0,0,0,0,{z_{{3}}}^{2},0,0,0,0,-z_{{6}}z_{{3}}+z_{{4}},2\,z_{{3}},0,0,0,0,0,0,0,0,0,0,{z_{{3}}}^{3},0,0,0,\\
&&-z_{{6}}{z_{{3}}}^{2}+3\,z_{{4}}z_{{3}},2\,{z_{{3}}}^{2},0,0,0,0,0,0,{z_{{3}}}^{4},0,0,-z_{{6}}{z_{{3}}}^{3}+6\,z_{{4}}{z_{{3}}}^{2},2\,{z_{{3}}}^{3},0,0,0,{z_{{3}}}^{5},0,-z_{{6}}{z_{{3}}}^{4}+10\,z_{{4}}{z_{{3}}}^{3},\\
&&2\,{z_{{3}}}^{4},0,{z_{{3}}}^{6},-z_{{6}}{z_{{3}}}^{5}+15\,z_{{4}}{z_{{3}}}^{4},2\,{z_{{3}}}^{5},21\,z_{{4}}{z_{{3}}}^{5},0,0,0]\\
X_{29}&=&[0,0,0,0,0,0,0,0,0,0,0,0,0,0,0,0,0,0,0,0,0,0,0,0,0,0,0,0,0,0,0,0,0,0,0,0,0,0,0,0,0,0,0,\\
&&0,0,0,0,0,0,0,0,0,0,0,0,0,0,0,1,0,0,0]\\
X_{30}&=&[0,0,0,0,0,0,0,0,0,0,0,0,0,0,0,0,0,0,0,0,0,0,0,0,0,0,0,0,0,0,0,0,0,0,0,0,0,0,0,0,0,0,0,\\
&&0,0,0,0,0,0,0,0,0,0,0,0,0,1,0,7\,z_{{3}},0,0,0]\\
X_{31}&=&[0,0,0,0,0,0,0,0,0,0,0,0,0,0,0,0,0,0,0,0,0,0,0,0,0,0,0,0,0,0,0,0,0,0,0,0,0,0,0,0,0,0,0,\\
&&0,0,0,0,0,0,0,0,0,1,0,0,0,6\,z_{{3}},0,21\,{z_{{3}}}^{2},0,0,0]\\
X_{32}&=&[0,0,0,0,0,0,0,0,0,0,0,0,0,0,0,0,0,0,0,0,0,0,0,0,0,0,0,0,0,0,0,0,0,0,0,0,0,0,0,0,0,0,0,\\
&&0,0,1,0,0,0,0,0,0,5\,z_{{3}},0,0,0,15\,{z_{{3}}}^{2},0,35\,{z_{{3}}}^{3},0,0,0]\\
X_{33}&=&[0,0,0,0,0,0,0,0,0,0,0,0,0,0,0,0,0,0,0,0,0,0,0,0,0,0,0,0,0,0,0,0,0,0,1,0,0,0,0,0,0,0,0,\\
&&0,0,4\,z_{{3}},0,0,0,0,0,0,10\,{z_{{3}}}^{2},0,0,0,20\,{z_{{3}}}^{3},0,35\,{z_{{3}}}^{4},0,0,0]\\
X_{34}&=&[0,0,0,0,0,0,0,0,0,0,0,0,0,0,0,0,0,0,1,0,0,0,0,0,0,0,0,0,0,0,0,0,0,0,3\,z_{{3}},0,0,0,0,\\
&&0,0,0,0,0,0,6\,{z_{{3}}}^{2},0,0,0,0,0,0,10\,{z_{{3}}}^{3},0,0,0,15\,{z_{{3}}}^{4},0,21\,{z_{{3}}}^{5},0,0,0]\\
X_{35}&=&[0,0,0,0,0,0,0,0,0,0,0,0,0,0,0,0,0,0,2\,z_{{3}},0,0,0,0,0,0,0,0,0,0,0,0,0,0,0,3\,{z_{{3}}}^{2},0,0,0,0,0,\\
&&0,0,0,0,0,4\,{z_{{3}}}^{3},0,0,0,0,0,0,5\,{z_{{3}}}^{4},0,0,0,6\,{z_{{3}}}^{5},0,7\,{z_{{3}}}^{6},0,0,0]\\
X_{36}&=&[0,0,0,0,0,0,0,0,0,0,0,0,0,0,0,0,0,0,{z_{{3}}}^{2},0,0,0,0,0,0,0,0,0,0,0,0,0,0,0,{z_{{3}}}^{3},0,0,0,0,0,\\
&&0,0,0,0,0,{z_{{3}}}^{4},0,0,0,0,0,0,{z_{{3}}}^{5},0,0,0,{z_{{3}}}^{6},0,{z_{{3}}}^{7},0,0,0]\\
\end{eqnarray*}
\end{footnotesize}
\begin{footnotesize}\tiny
\begin{eqnarray*}
T_{22}&=&[0,0,0,0,0,0,0,0,0,0,0,0,0,0,0,0,0,0,0,0,0,0,0,0,0,0,0,0,0,0,0,0,0,0,0,0,0,0,0,0,0,0,0,0,\\
&&0,0,0,0,0,0,0,0,0,0,0,-z_{{3}},0,-1,-7\,z_{{4}},0,0,0]\\
T_{23}&=&[0,0,0,0,0,0,0,0,0,0,0,0,0,0,0,0,0,0,0,0,0,0,0,0,0,0,0,0,0,0,0,0,0,0,0,0,0,0,0,0,0,0,0,0,\\
&&0,0,0,0,0,0,-z_{{3}},0,0,-1,0,-6\,{z_{{3}}}^{2},-6\,z_{{4}}+z_{{6}}z_{{3}},-8\,z_{{3}},-63\,z_{{4}}z_{{3}},0,0,0]\\
T_{24}&=&[0,0,0,0,0,0,0,0,0,0,0,0,0,0,0,0,0,0,0,0,0,0,0,0,0,0,0,0,0,0,0,0,0,0,0,0,0,0,0,0,0,0,-z_{{3}},\\
&&0,0,0,-1,0,0,0,-5\,{z_{{3}}}^{2},0,-5\,z_{{4}}+z_{{6}}z_{{3}},-7\,z_{{3}},0,-15\,{z_{{3}}}^{3},-45\,z_{{4}}z_{{3}}+5\,z_{{6}}{z_{{3}}}^{2},-25\,{z_{{3}}}^{2},-210\,z_{{4}}{z_{{3}}}^{2},0,0,0]\\
T_{25}&=&[0,0,0,0,0,0,0,0,0,0,0,0,0,0,0,0,0,0,0,0,0,0,0,0,0,0,0,0,0,0,-z_{{3}},0,0,0,0,-1,0,0,0,0,0,0,-4\,{z_{{3}}}^{2},\\
&&0,0,z_{{6}}z_{{3}}-4\,z_{{4}},-6\,z_{{3}},0,0,0,-10\,{z_{{3}}}^{3},0,4\,z_{{6}}{z_{{3}}}^{2}-30\,z_{{4}}z_{{3}},-18\,{z_{{3}}}^{2},0,-20\,{z_{{3}}}^{4},-120\,z_{{4}}{z_{{3}}}^{2}+10\,z_{{6}}{z_{{3}}}^{3},\\
&&-40\,{z_{{3}}}^{3},-350\,z_{{4}}{z_{{3}}}^{3},0,0,0]\\
T_{26}&=&[0,0,0,0,0,0,0,0,0,0,0,0,0,-z_{{3}},0,0,0,0,0,-1,0,0,0,0,0,0,0,0,0,0,-3\,{z_{{3}}}^{2},0,0,0,-3\,z_{{4}}+z_{{6}}z_{{3}},-5\,z_{{3}},\\
&&0,0,0,0,0,0,-6\,{z_{{3}}}^{3},0,0,-18\,z_{{4}}z_{{3}}+3\,z_{{6}}{z_{{3}}}^{2},-12\,{z_{{3}}}^{2},0,0,0,-10\,{z_{{3}}}^{4},0,6\,z_{{6}}{z_{{3}}}^{3}-60\,z_{{4}}{z_{{3}}}^{2},-22\,{z_{{3}}}^{3},0,\\
&&-15\,{z_{{3}}}^{5},10\,z_{{6}}{z_{{3}}}^{4}-150\,z_{{4}}{z_{{3}}}^{3},-35\,{z_{{3}}}^{4},-315\,z_{{4}}{z_{{3}}}^{4},0,0,0]\\
T_{27}&=&[0,0,0,0,0,0,0,0,0,0,0,0,0,-2\,{z_{{3}}}^{2},0,0,0,0,-2\,z_{{4}}+z_{{6}}z_{{3}},-4\,z_{{3}},0,0,0,0,0,0,0,0,0,0,-3\,{z_{{3}}}^{3},0,0,0,\\
&&2\,z_{{6}}{z_{{3}}}^{2}-9\,z_{{4}}z_{{3}},-7\,{z_{{3}}}^{2},0,0,0,0,0,0,-4\,{z_{{3}}}^{4},0,0,-24\,z_{{4}}{z_{{3}}}^{2}+3\,z_{{6}}{z_{{3}}}^{3},-10\,{z_{{3}}}^{3},0,0,0,-5\,{z_{{3}}}^{5},0,\\
&&4\,z_{{6}}{z_{{3}}}^{4}-50\,z_{{4}}{z_{{3}}}^{3},-13\,{z_{{3}}}^{4},0,-6\,{z_{{3}}}^{6},5\,z_{{6}}{z_{{3}}}^{5}-90\,z_{{4}}{z_{{3}}}^{4},-16\,{z_{{3}}}^{5},-147\,z_{{4}}{z_{{3}}}^{5},0,0,0]\\
T_{28}&=&[0,0,0,0,0,0,0,0,0,0,0,0,0,-{z_{{3}}}^{3},0,0,0,0,z_{{6}}{z_{{3}}}^{2}-3\,z_{{4}}z_{{3}},-3\,{z_{{3}}}^{2},0,0,0,0,0,0,0,0,0,0,-{z_{{3}}}^{4},0,0,0,\\
&&z_{{6}}{z_{{3}}}^{3}-6\,z_{{4}}{z_{{3}}}^{2},-3\,{z_{{3}}}^{3},0,0,0,0,0,0,-{z_{{3}}}^{5},0,0,-10\,z_{{4}}{z_{{3}}}^{3}+z_{{6}}{z_{{3}}}^{4},-3\,{z_{{3}}}^{4},0,0,0,-{z_{{3}}}^{6},0,\\
&&-15\,z_{{4}}{z_{{3}}}^{4}+z_{{6}}{z_{{3}}}^{5},-3\,{z_{{3}}}^{5},0,-{z_{{3}}}^{7},-21\,z_{{4}}{z_{{3}}}^{5}+z_{{6}}{z_{{3}}}^{6},-3\,{z_{{3}}}^{6},-28\,z_{{4}}{z_{{3}}}^{6},0,0,0]\\
T_{29}&=&[0,0,0,0,0,0,0,0,0,0,0,0,0,0,0,0,0,0,0,0,0,0,0,0,0,0,0,0,0,0,0,0,0,0,0,0,0,0,0,0,0,0,0,0,\\
&&0,0,0,0,0,0,0,0,0,0,0,0,0,0,-z_{{3}},0,0,0]\\
T_{30}&=&[0,0,0,0,0,0,0,0,0,0,0,0,0,0,0,0,0,0,0,0,0,0,0,0,0,0,0,0,0,0,0,0,0,0,0,0,0,0,0,0,0,0,0,0,\\
&&0,0,0,0,0,0,0,0,0,0,0,0,-z_{{3}},0,-7\,{z_{{3}}}^{2},0,0,0]\\
T_{31}&=&[0,0,0,0,0,0,0,0,0,0,0,0,0,0,0,0,0,0,0,0,0,0,0,0,0,0,0,0,0,0,0,0,0,0,0,0,0,0,0,0,0,0,0,0,\\
&&0,0,0,0,0,0,0,0,-z_{{3}},0,0,0,-6\,{z_{{3}}}^{2},0,-21\,{z_{{3}}}^{3},0,0,0]\\
T_{32}&=&[0,0,0,0,0,0,0,0,0,0,0,0,0,0,0,0,0,0,0,0,0,0,0,0,0,0,0,0,0,0,0,0,0,0,0,0,0,0,0,0,0,0,0,0,\\
&&0,-z_{{3}},0,0,0,0,0,0,-5\,{z_{{3}}}^{2},0,0,0,-15\,{z_{{3}}}^{3},0,-35\,{z_{{3}}}^{4},0,0,0]\\
T_{33}&=&[0,0,0,0,0,0,0,0,0,0,0,0,0,0,0,0,0,0,0,0,0,0,0,0,0,0,0,0,0,0,0,0,0,0,-z_{{3}},0,0,0,0,0,0,0,0\\
&&,0,0,-4\,{z_{{3}}}^{2},0,0,0,0,0,0,-10\,{z_{{3}}}^{3},0,0,0,-20\,{z_{{3}}}^{4},0,-35\,{z_{{3}}}^{5},0,0,0]\\
T_{34}&=&[0,0,0,0,0,0,0,0,0,0,0,0,0,0,0,0,0,0,-z_{{3}},0,0,0,0,0,0,0,0,0,0,0,0,0,0,0,-3\,{z_{{3}}}^{2},0,0,0,0,0,\\
&&0,0,0,0,0,-6\,{z_{{3}}}^{3},0,0,0,0,0,0,-10\,{z_{{3}}}^{4},0,0,0,-15\,{z_{{3}}}^{5},0,-21\,{z_{{3}}}^{6},0,0,0]\\
T_{35}&=&[0,0,0,0,0,0,0,0,0,0,0,0,0,0,0,0,0,0,-2\,{z_{{3}}}^{2},0,0,0,0,0,0,0,0,0,0,0,0,0,0,0,-3\,{z_{{3}}}^{3},0,0,0,0\\
&&,0,0,0,0,0,0,-4\,{z_{{3}}}^{4},0,0,0,0,0,0,-5\,{z_{{3}}}^{5},0,0,0,-6\,{z_{{3}}}^{6},0,-7\,{z_{{3}}}^{7},0,0,0]\\
T_{36}&=&[0,0,0,0,0,0,0,0,0,0,0,0,0,0,0,0,0,0,-{z_{{3}}}^{3},0,0,0,0,0,0,0,0,0,0,0,0,0,0,0,-{z_{{3}}}^{4},0,0,0,0,0,\\
&&0,0,0,0,0,-{z_{{3}}}^{5},0,0,0,0,0,0,-{z_{{3}}}^{6},0,0,0,-{z_{{3}}}^{7},0,-{z_{{3}}}^{8},0,0,0]\\
\end{eqnarray*}
\end{footnotesize}
\section*{Appendix B: The differential operators of the homogeneous linear system of
PDEs (\ref{c13}) }
\begin{small}
\begin{equation*}
\begin{array}{ll}
X_{16}=&[0,0,0,0,0,0,0,0,0,0,0,0,0,0,0,0,0,0,0,0,0,0,0,0,0,0,0,0,0,0,0,0,0,\\
&0,0,0,0,0,0,0,1,0,0,0,0,0]\\
X_{17}=&[0,0,0,0,0,0,0,0,0,0,0,0,0,0,0,0,0,0,0,0,0,0,0,0,0,0,0,0,0,0,0,0,0,\\
&0,1,0,0,0,0,0,5\,l_{{3}},2,0,0,0,0]\\
X_{18}=&[0,0,0,0,0,0,0,0,0,0,0,0,0,0,0,0,0,0,0,0,0,0,0,0,1,0,0,0,0,0,0,0,0,\\
&0,4\,l_{{3}},0,2,0,0,0,10\,{l_{{3}}}^{2},8\,l_{{3}},2,0,0,0]\\
X_{19}=&[0,0,0,0,0,0,0,0,0,1,0,0,0,0,0,0,0,0,0,0,0,0,0,0,3\,l_{{3}},0,0,2,0,0,0,0,\\
&0,0,6\,{l_{{3}}}^{2},0,6\,l_{{3}},0,2,0,10\,{l_{{3}}}^{3},12\,{l_{{3}}}^{2},6\,l_{{3}},0,0,0]\\
X_{20}=&[0,0,0,0,0,0,0,0,0,2\,l_{{3}},0,0,0,2,0,0,0,0,0,0,0,0,0,0,3\,{l_{{3}}}^{2},0,0,4\,l_{{3}},0,0,\\
&2,0,0,0,4\,{l_{{3}}}^{3},0,6\,{l_{{3}}}^{2},0,4\,l_{{3}},0,5\,{l_{{3}}}^{4},8\,{l_{{3}}}^{3},6\,{l_{{3}}}^{2},0,0,0]\\
X_{21}=&[0,0,0,0,0,0,0,0,0,{l_{{3}}}^{2},0,0,0,2\,l_{{3}},0,0,0,2,0,0,0,0,0,0,{l_{{3}}}^{3},0,0,2\,{l_{{3}}}^{2},0,0,\\
&2\,l_{{3}},0,0,0,{l_{{3}}}^{4},0,2\,{l_{{3}}}^{3},0,2\,{l_{{3}}}^{2},0,{l_{{3}}}^{5},2\,{l_{{3}}}^{4},2\,{l_{{3}}}^{3},0,0,0]\\
T_{16}=&[0,0,0,0,0,0,0,0,0,0,0,0,0,0,0,0,0,0,0,0,0,0,0,0,0,0,0,0,0,0,0,0,0,\\
&0,0,0,0,0,0,0,-l_{{3}},-1,0,0,0,0]\\
T_{17}=&[0,0,0,0,0,0,0,0,0,0,0,0,0,0,0,0,0,0,0,0,0,0,0,0,0,0,0,0,0,0,0,0,0,\\
&0,-l_{{3}},0,-1,0,0,0,-5\,{l_{{3}}}^{2},-7\,l_{{3}},-4,0,0,0]\\
T_{18}=&[0,0,0,0,0,0,0,0,0,0,0,0,0,0,0,0,0,0,0,0,0,0,0,0,-l_{{3}},0,0,-1,0,0,0,\\
&0,0,0,-4\,{l_{{3}}}^{2},0,-6\,l_{{3}},0,-4,0,-10\,{l_{{3}}}^{3},-18\,{l_{{3}}}^{2},-18\,l_{{3}},0,0,0]\\
T_{19}=&[0,0,0,0,0,0,0,0,0,-l_{{3}},0,0,0,-1,0,0,0,0,0,0,0,0,0,0,-3\,{l_{{3}}}^{2},0,0,-5\,l_{{3}},\\
&0,0,-4,0,0,0,-6\,{l_{{3}}}^{3},0,-12\,{l_{{3}}}^{2},0,-14\,l_{{3}},0,-10\,{l_{{3}}}^{4},-22\,{l_{{3}}}^{3},-30\,{l_{{3}}}^{2},0,0,0]\\
T_{20}=&[0,0,0,0,0,0,0,0,0,-2\,{l_{{3}}}^{2},0,0,0,-4\,l_{{3}},0,0,0,-4,0,0,0,0,0,0,-3\,{l_{{3}}}^{3},0,0,\\
&-7\,{l_{{3}}}^{2},0,0,-10\,l_{{3}},0,0,0,-4\,{l_{{3}}}^{4},0,-10\,{l_{{3}}}^{3},0,-16\,{l_{{3}}}^{2},0,-5\,{l_{{3}}}^{5},-13\,{l_{{3}}}^{4},-22\,{l_{{3}}}^{3},0,0,0]\\
T_{21}=&[0,0,0,0,0,0,0,0,0,-{l_{{3}}}^{3},0,0,0,-3\,{l_{{3}}}^{2},0,0,0,-6\,l_{{3}},0,0,0,0,0,0,-{l_{{3}}}^{4},0,0,\\
&-3\,{l_{{3}}}^{3},0,0,-6\,{l_{{3}}}^{2},0,0,0,-{l_{{3}}}^{5},0,-3\,{l_{{3}}}^{4},0,-6\,{l_{{3}}}^{3},0,-{l_{{3}}}^{6},-3\,{l_{{3}}}^{5},-6\,{l_{{3}}}^{4},0,0,0]\\
\end{array}
\end{equation*}
\end{small}
\section*{Appendix C: The differential operators of the homogeneous linear system of
PDEs (\ref{c15}) }
\begin{small}
\begin{equation*}
\begin{array}{ll}
X_{11}=&[0,0,0,0,0,0,0,0,0,0,0,0,0,0,0,0,0,0,0,0,0,0,0,0,0,0,0,1,0,0,0,0,0,0,0,0]\\
X_{12}=&[0,0,0,0,0,0,0,0,0,0,0,0,0,0,0,0,0,0,0,1,0,0,0,0,0,0,0,4\,m_{{3}},2,0,-m_{{12}},0,\\
&3\,m_{{9}},0,0,0]\\
X_{13}=&[0,0,0,0,0,0,1,0,0,0,0,0,0,0,0,0,0,0,0,3\,m_{{3}},0,2,0,0,0,0,0,6\,{m_{{3}}}^{2},6\,m_{{3}},2,\\
&3\,m_{{9}}-3\,m_{{3}}m_{{12}},0,-6\,m_{{6}}+9\,m_{{3}}m_{{9}},0,0,0]\\
X_{14}=&[0,0,0,0,0,0,2\,m_{{3}},0,0,2,0,0,0,0,0,0,0,0,0,3\,{m_{{3}}}^{2},0,4\,m_{{3}},0,2,0,0,0,4\,{m_{{3}}}^{3},\\
&6\,{m_{{3}}}^{2},4\,m_{{3}},-6\,m_{{6}}+6\,m_{{3}}m_{{9}}-3\,{m_{{3}}}^{2}m_{{12}},0,-12\,m_{{6}}m_{{3}}+9\,{m_{{3}}}^{2}m_{{9}},0,0,0]\\
X_{15}=&[0,0,0,0,0,0,{m_{{3}}}^{2},0,0,2\,m_{{3}},0,0,2,0,0,0,0,0,0,{m_{{3}}}^{3},0,2\,{m_{{3}}}^{2},0,2\,m_{{3}},0,0,0,{m_{{3}}}^{4},\\
&2\,{m_{{3}}}^{3},2\,{m_{{3}}}^{2},-6\,m_{{6}}m_{{3}}+3\,{m_{{3}}}^{2}m_{{9}}-{m_{{3}}}^{3}m_{{12}},0,-6\,m_{{6}}{m_{{3}}}^{2}+3\,{m_{{3}}}^{3}m_{{9}},0,0,0]\\
T_{11}=&[0,0,0,0,0,0,0,0,0,0,0,0,0,0,0,0,0,0,0,0,0,0,0,0,0,0,0,-m_{{3}},-1,0,0,0,0,0,0,0]\\
T_{12}=&[0,0,0,0,0,0,0,0,0,0,0,0,0,0,0,0,0,0,0,-m_{{3}},0,-1,0,0,0,0,0,-4\,{m_{{3}}}^{2},-6\,m_{{3}},-4,\\
&m_{{3}}m_{{12}},0,-3\,m_{{3}}m_{{9}},0,0,0]\\
T_{13}=&[0,0,0,0,0,0,-m_{{3}},0,0,-1,0,0,0,0,0,0,0,0,0,-3\,{m_{{3}}}^{2},0,-5\,m_{{3}},0,-4,0,0,0,-6\,{m_{{3}}}^{3},\\
&-12\,{m_{{3}}}^{2},-14\,m_{{3}},-3\,m_{{3}}m_{{9}}+3\,{m_{{3}}}^{2}m_{{12}},-6,6\,m_{{6}}m_{{3}}-9\,{m_{{3}}}^{2}m_{{9}},0,0,0]\\
T_{14}=&[0,0,0,0,0,0,-2\,{m_{{3}}}^{2},0,0,-4\,m_{{3}},0,0,-4,0,0,0,0,0,0,-3\,{m_{{3}}}^{3},0,-7\,{m_{{3}}}^{2},0,-10\,m_{{3}},\\
&0,-6,0,-4\,{m_{{3}}}^{4},-10\,{m_{{3}}}^{3},-16\,{m_{{3}}}^{2},6\,m_{{6}}m_{{3}}-6\,{m_{{3}}}^{2}m_{{9}}+3\,{m_{{3}}}^{3}m_{{12}},-12\,m_{{3}},\\
&12\,m_{{6}}{m_{{3}}}^{2}-9\,{m_{{3}}}^{3}m_{{9}},0,0,0]\\
T_{15}=&[0,0,0,0,0,0,-{m_{{3}}}^{3},0,0,-3\,{m_{{3}}}^{2},0,0,-6\,m_{{3}},0,0,-6,0,0,0,-{m_{{3}}}^{4},0,-3\,{m_{{3}}}^{3},0,-6\,{m_{{3}}}^{2},\\
&0,-6\,m_{{3}},0,-{m_{{3}}}^{5},-3\,{m_{{3}}}^{4},-6\,{m_{{3}}}^{3},6\,m_{{6}}{m_{{3}}}^{2}-3\,{m_{{3}}}^{3}m_{{9}}+{m_{{3}}}^{4}m_{{12}},-6\,{m_{{3}}}^{2},\\
&6\,m_{{6}}{m_{{3}}}^{3}-3\,{m_{{3}}}^{4}m_{{9}},0,0,0]\\
\end{array}
\end{equation*}
\end{small}
\section*{Appendix D: The differential operators of the homogeneous linear system of
PDEs (\ref{c17}) }
\begin{small}
\begin{equation*}
\begin{array}{ll}
X_{7}=&[0,0,0,0,0,0,0,0,0,0,0,0,0,0,1,0,0,0,0,0,0,0,0,0,0,0]\\
X_{8}=&[0,0,0,0,1,0,0,0,0,0,0,0,0,0,3\,n_{{3}},2,0,0,0,3\,n_{{8}},-2\,n_{{18}}+6\,n_{{9}},0,5\,n_{{17}}-8\,n_{{7}}+3\,n_{{8}}n_{{6}},\\
&0,0,0]\\
X_{9}=&[0,0,0,0,2\,n_{{3}},0,2,0,0,0,0,0,0,0,3\,{n_{{3}}}^{2},4\,n_{{3}},2,0,0,-6\,n_{{6}}+6\,n_{{3}}n_{{8}},12\,n_{{3}}n_{{9}}-4\,n_{{3}}n_{{18}}\\
&-14\,n_{{7}}+2\,n_{{17}},0,-6\,{n_{{6}}}^{2}+6\,n_{{8}}n_{{6}}n_{{3}}-16\,n_{{3}}n_{{7}}+10\,n_{{3}}n_{{17}}+6\,n_{{5}}-6\,n_{{16}},0,0,0]\\
X_{10}=&[0,0,0,0,{n_{{3}}}^{2},0,2\,n_{{3}},0,2,0,0,0,0,0,{n_{{3}}}^{3},2\,{n_{{3}}}^{2},2\,n_{{3}},0,0,-6\,n_{{6}}n_{{3}}+3\,{n_{{3}}}^{2}n_{{8}},-2\,{n_{{3}}}^{2}n_{{18}}\\
&+6\,{n_{{3}}}^{2}n_{{9}}+2\,n_{{3}}n_{{17}}-14\,n_{{3}}n_{{7}}+12\,n_{{5}},0,-6\,{n_{{6}}}^{2}n_{{3}}+3\,n_{{8}}n_{{6}}{n_{{3}}}^{2}+5\,{n_{{3}}}^{2}n_{{17}}-8\,{n_{{3}}}^{2}n_{{7}}\\
&+6\,n_{{3}}n_{{5}}-6\,n_{{3}}n_{{16}},0,0,0]\\
T_{7}=&[0,0,0,0,0,0,0,0,0,0,0,0,0,0,-n_{{3}},-1,0,0,0,0,0,n_{{10}},4\,n_{{7}}-n_{{17}},0,0,0]\\
T_{8}=&[0,0,0,0,-n_{{3}},0,-1,0,0,0,0,0,0,0,-3\,{n_{{3}}}^{2},-5\,n_{{3}},-4,0,0,-3\,n_{{3}}n_{{8}},-6\,n_{{3}}n_{{9}}+2\,n_{{3}}n_{{18}}\\
&+4\,n_{{7}}-n_{{17}},-3\,n_{{8}}+3\,n_{{3}}n_{{10}},-3\,n_{{8}}n_{{6}}n_{{3}}+20\,n_{{3}}n_{{7}}-8\,n_{{3}}n_{{17}}-12\,n_{{5}},0,0,0]\\
T_{9}=&[0,0,0,0,-2\,{n_{{3}}}^{2},0,-4\,n_{{3}},0,-4,0,0,0,0,0,-3\,{n_{{3}}}^{3},-7\,{n_{{3}}}^{2},-10\,n_{{3}},-6,0,6\,n_{{6}}n_{{3}}\\
&-6\,{n_{{3}}}^{2}n_{{8}},-12\,{n_{{3}}}^{2}n_{{9}}+4\,{n_{{3}}}^{2}n_{{18}}-4\,n_{{3}}n_{{17}}+22\,n_{{3}}n_{{7}}-12\,n_{{5}},6\,n_{{6}}-6\,n_{{3}}n_{{8}}+3\,{n_{{3}}}^{2}n_{{10}},\\
&6\,{n_{{6}}}^{2}n_{{3}}-6\,n_{{8}}n_{{6}}{n_{{3}}}^{2}-13\,{n_{{3}}}^{2}n_{{17}}+28\,{n_{{3}}}^{2}n_{{7}}-30\,n_{{3}}n_{{5}}+6\,n_{{3}}n_{{16}},0,0,0]\\
T_{10}=&[0,0,0,0,-{n_{{3}}}^{3},0,-3\,{n_{{3}}}^{2},0,-6\,n_{{3}},0,-6,0,0,0,-{n_{{3}}}^{4},-3\,{n_{{3}}}^{3},-6\,{n_{{3}}}^{2},-6\,n_{{3}},0,6\,n_{{6}}{n_{{3}}}^{2}\\
&-3\,{n_{{3}}}^{3}n_{{8}},-6\,{n_{{3}}}^{3}n_{{9}}+2\,{n_{{3}}}^{3}n_{{18}}+18\,{n_{{3}}}^{2}n_{{7}}-3\,{n_{{3}}}^{2}n_{{17}}-24\,n_{{3}}n_{{5}},-3\,{n_{{3}}}^{2}n_{{8}}+6\,n_{{6}}n_{{3}}\\
&+{n_{{3}}}^{3}n_{{10}},6\,{n_{{6}}}^{2}{n_{{3}}}^{2}-3\,{n_{{3}}}^{3}n_{{8}}n_{{6}}+12\,{n_{{3}}}^{3}n_{{7}}+6\,{n_{{3}}}^{2}n_{{16}}-18\,{n_{{3}}}^{2}n_{{5}}-6\,{n_{{3}}}^{3}n_{{17}},0,0,0]\\
\end{array}
\end{equation*}
\end{small}
\section*{Appendix E: The differential operators of the homogeneous linear system of
PDEs (\ref{c19}) }
\begin{small}
\begin{equation*}
\begin{array}{ll}
X_{1}=&[0,1,0,0,0,0,0,0,0,0,0,0,0,0,0,0,0,0]\\
X_{2}=&[0,0,1,0,0,0,0,0,0,0,0,0,0,0,0,0,t_{{16}},0]\\
X_{3}=&[0,0,t_{{3}},t_{{4}},0,-t_{{6}},-2\,t_{{7}},-3\,t_{{8}},-4\,t_{{9}},-4\,t_{{10}},-3\,t_{{11}},-t_{{12}},-2\,t_{{13}},-2\,t_{{14}},-t_{{15}},0,t_{{17}},t_{{18}}]\\
X_{4}=&[0,0,0,1,0,0,0,0,0,0,0,0,0,0,0,0,0,t_{{16}}]\\
X_{5}=&[0,0,0,2\,t_{{3}},2,0,0,0,-t_{{10}},0,-3\,t_{{8}},0,-t_{{14}},0,t_{{12}},0,0,t_{{17}}+t_{{16}}t_{{3}}]\\
X_{6}=&[0,0,0,{t_{{3}}}^{2},2\,t_{{3}},2,0,0,-3\,t_{{8}}-t_{{3}}t_{{10}},0,-3\,t_{{3}}t_{{8}},0,-t_{{3}}t_{{14}}-t_{{12}},0,t_{{3}}t_{{12}},0,0,t_{{17}}t_{{3}}]\\
T_{1}=&[1,0,0,0,0,0,0,0,0,0,0,0,0,0,0,0,0,0]\\
T_{2}=&[0,0,-t_{{3}},-2\,t_{{4}},-t_{{5}},0,t_{{7}},2\,t_{{8}},2\,t_{{9}},3\,t_{{10}},t_{{11}},-2\,t_{{12}},-2\,t_{{13}},-t_{{14}},-3\,t_{{15}},t_{{16}},0,-t_{{18}}]\\
T_{3}=&[0,0,-{t_{{3}}}^{2},-3\,t_{{4}}t_{{3}},-t_{{5}}t_{{3}}-3\,t_{{4}},-4\,t_{{5}}+t_{{3}}t_{{6}},-3\,t_{{6}}+3\,t_{{3}}t_{{7}},5\,t_{{3}}t_{{8}},-t_{{11}}+t_{{4}}t_{{10}}+6\,t_{{3}}t_{{9}},\\
&5\,t_{{8}}+7\,t_{{3}}t_{{10}},5\,t_{{4}}t_{{8}}+4\,t_{{3}}t_{{11}},-t_{{3}}t_{{12}},t_{{5}}t_{{12}}+t_{{4}}t_{{14}}-t_{{15}},t_{{3}}t_{{14}}-t_{{12}},-4\,t_{{4}}t_{{12}}-2\,t_{{3}}t_{{15}},t_{{17}},0,-2\,t_{{18}}t_{{3}}]\\
T_{4}=&[0,0,0,-t_{{3}},-1,0,0,0,0,0,2\,t_{{8}},0,0,0,-3\,t_{{12}},0,0,-t_{{16}}t_{{3}}]\\
T_{5}=&[0,0,0,-2\,{t_{{3}}}^{2},-4\,t_{{3}},-4,0,0,t_{{3}}t_{{10}}+2\,t_{{8}},0,7\,t_{{3}}t_{{8}},0,t_{{3}}t_{{14}}-2\,t_{{12}},0,-7\,t_{{3}}t_{{12}},0,0,-t_{{16}}{t_{{3}}}^{2}-t_{{17}}t_{{3}}]\\
T_{6}=&[0,0,0,-{t_{{3}}}^{3},-3\,{t_{{3}}}^{2},-6\,t_{{3}},-6,0,{t_{{3}}}^{2}t_{{10}}+5\,t_{{3}}t_{{8}},0,5\,{t_{{3}}}^{2}t_{{8}},0,-t_{{3}}t_{{12}}+{t_{{3}}}^{2}t_{{14}},0,-4\,{t_{{3}}}^{2}t_{{12}},0,0,-t_{{17}}{t_{{3}}}^{2}]\\
\end{array}
\end{equation*}
\end{small}

\end{document}